\documentclass[12pt]{amsart}
\usepackage{color}
\usepackage{epsfig}
\usepackage{verbatim}

\begin{document}
\title {Growth rate of binary matroids with no $P_9^*$-minor}

\maketitle 

\begin {center}
S. R. Kingan 
\footnote{The author is partially supported by  PSC-CUNY grant number 66305-00 45.}\\     
Department of Mathematics \\
Brooklyn College, City University of New York\\
 Brooklyn, NY 11210\\
skingan@brooklyn.cuny.edu\\  
\end {center}
\bigskip

\begin{abstract}  We prove that the non-regular binary matroids with no $P_9^*$-minor have linear growth rate and the maximum size binary matroids with no $P_9^*$-minor are graphic. The main technique in the proof is the Strong Splitter Theorem using which we find the precise infinite families of 3-connected binary matroids with no $P_9^*$-minor.
\end{abstract}

\bigskip 

\section {\bf Introduction}

 The {\it growth rate function} of a minor-closed class $\mathcal M$, denoted by $h_{\mathcal M}(r)$, is defined as the maximum number of elements in a simple 3-connected rank-$r$ matroid in $\mathcal M$ if the number is finite and infinity otherwise. We put 3-connected in the definition itself since 
matroids that are not 3-connected may be built up from 3-connected matroids using direct-sums and 2-sums [\ref{Oxley1987}, 8.3.1].  
For example, if $\mathcal M$ is the class of graphic matroids, then $\mathcal M$ has just one rank-$r$ infinite family of matroids, {\it namely}, the complete graphs of rank-$r$ denoted by $K_{r+1}$, for $r\ge 1$. The growth rate function is the size of  $K_{r+1}$, which is $\frac {r(r+1)}{2}$. We call $K_{r+1}$ an infinite family of rank-$r$ extremal matroids in $\mathcal M$.
 Similarly, if $\mathcal M$ is the class of binary matroids, then again $\mathcal M$ has just one rank-$r$ infinite family of extremal matroids, {\it namely}, the  rank-$r$ projective geometry $PG(r-1, 2)$, for $r\ge 2$. Note that the rank-$r$ complete graph $K_{r+1}$ is a minor of the rank-$r$ projective geometry $PG(r-1, 2)$. The growth rate function is the size of  $PG(r-1, 2)$, which is $2^r-1$.  

A non-trivial example would be Oxley's characterization of  the  binary matroids with no $P_9$ or $P_9^*$-minor [\ref{Oxley1987}]. This class has one  rank-$r$ infinite family of matroids, {\it namely}, the rank-$r$ binary spikes $Z_r$, for $r\ge 4$. The matroid $Z_r$ can be represented by the matrix $[I_r | D]$, where $D$ has $r+1$  columns with zeros on the diagonal and ones elsewhere. The growth rate function of this class is the size of  $Z_r$, which is $(2r+1)$. 

A minor-closed class may have more than one infinite family of 3-connected extremal matroids.  For example, if $\mathcal M$ is the class of graphs with no minor isomorphic to the prism graph $(K_5\backslash e)^*$, then Dirac proved that  the two distinct recursively defined rank-$r$ infinite families of 3-connected extremal graphs are $W_r$, for $r\ge 3$, and $K_{3, p}'''$, for $p\ge 3$ [\ref{Dirac1963}]. Observe that $W_r$ and $K_{3, p}'''$ are not deletion-minors of each other as in the case of $K_{r+1}$ and $PG(r-1, 2)$. The size of $W_r$ is $2r+1$ and the size of $K_{3, p}'''$ is $3r-3$. Since the  growth-rate function  is the size of the largest  rank-$r$ matroid, for this class it is $3r-3$. 
Another example appears in [\ref{KinganLemos2002}, 8.2]. The class of $3$-connected binary almost-regular matroid with no $E_5$-minor has three distinct infinite extremal families of matroids $\mathcal F_1(m, n, r)$ and $\mathcal F_2(m, n, r)$, for $m, n, r\ge 1$, and $S_{3n+1}$, for $n\ge 3$. 

The above examples illustrate the need for a careful definition of extremal matroids. Let $\mathcal M$ be a minor-closed class of matroids. Suppose $\mathcal M$ has $k\ge 1$ infinite families of recursively defined  rank-$r$ 3-connected matroids $M^1_r, \dots M^k_r$, where $M^i_r$ is not a deletion-minor of $M^j_r$ for $1\le i, j \le k$. Then, we say the class has  {\it $k$ distinct  extremal families of matroids}. Each extremal family has a growth-rate function $f(r)$ and the growth-rate function for the class $\mathcal M$, $h_{\mathcal M}(r)$, is the largest $f(r)$.

In this paper we give a complete characterization of the binary matroids with no $P_9^*$-minor. Matrix representations for $P_9$ and its dual $P_9^*$ are given below.   
\[ 
P_{9}=\left[ 
\begin{array}{c|ccccc}
&   0&1&1&1&1 \\
I_4&1&0&1&1&1 \\
&   1&1&0&1&0 \\
&   1&1&1&1&0
\end{array} 
\right] 
P_9^*=\left[ 
\begin{array}{c|cccc}
&    0&1&1&1 \\
&    1&0&1&1 \\
I_5& 1&1&0&1 \\
&    1&1&1&1 \\
&    1&1&0&0
\end{array} 
\right] 
\]
\noindent In addition to the binary spikes this class has the infinite family of matroids  $\Omega _r$, for $r\ge 5$, which has size $4r-5$. Whereas the technique used in Oxley's paper is the Splitter Theorem, the result in this paper requires the Strong Splitter Theorem [\ref{KinganLemos2014}]. 
The next theorem is the main theorem in this paper.
\bigskip

\noindent {\bf Theorem 1.1.} Let $M$ be a simple rank-$r$ binary matroid with no $P_9^*$ minor. Then,  $|E(M)| \le \frac {r(r+1)}{2}$ with this bound being attained if and only if $M\cong K_{r+1}$. Moreover, the growth rate function of non-regular matroids in $EX[P_9^*]$ is $4r-5$, for $r\ge 6$. $\qed$
\bigskip

The above result may be added to the limited set of results of its type. Notably, Kung, Mayhew, Pivotto, and Royle proved that $EX[AG(3,2)]$ has growth rate function  $\frac {r(r+1)}{2}$ [\ref{KungMayhewPivottoRoylesubmitted}], however, their proof has a computer component to it. It is important to note that,  unlike their proof, Theorem 1.1 has no computer component. 
\footnote {Note to referee: There is a lot of confusion in the matroid community about use of software. Papers that use Macek have to cite a computer component because of the way Macek works - Macek checks whether or not a matroid has a minor but does not give the explicit sequence of deletions and contractions to get the minor. When the explicit deletions and contractions are found there is no ``computer component." Moreover, those who use Macek are brute-force checking hundreds of matroids and it is not practical to do anything by hand. Hence they have to verify their computations with other software, which is the standard in computer science (a result is acceptable in computer science if two different programs give the same result). None of this is relevant to my work. Most of the few tables in the Appendix appear in my dissertation where linear transformations were found by hand between isomorphic binary matroids.} Moreover, the technique is completely different.  

Mader proved that the growth rate function of any proper minor-closed class of graphs is a linear function of the rank $r$ [\ref{Oxley2012}, 14.10.1].    Geelen, Kung and Whittle  proved that if a proper class of binary matroids $\mathcal M$ contains the entire class of graphs, then the growth rate function of $\mathcal M$ is quadratic [\ref{Oxley2012}, 14.10.7]. 
Our method for finding the infinite families of non-regular extremal matroids in $EX[P_9^*]$ illustrates to some extent why the growth rates of proper minor closed classes of binary matroids exhibit the bifurcation between linear and quadratic.

\section {The size of $\Omega_r$}

Before presenting the infinite family $\Omega _r$ let us make a new definition that gives us language for describing $\Omega _r$ based on the Strong Splitter Theorem [\ref{KinganLemos2014}, Theorem 1.4] given below: . 

\bigskip
\noindent{\bf Theorem 2.1.} {\it (Strong Splitter Theorem)
Suppose $N$ is a $3$-connected proper minor of a $3$-connected matroid $M$ such that,  if  $N$ is a wheel or a whirl then $M$ has no larger wheel or whirl-minor, respectively. Further, suppose $m=r(M)-r(N)$. Then there is a sequence of $3$-connected matroids $M_0,M_1,\dots,M_n$, for some integer $n\ge m$, such that:
\begin{enumerate}
\item[(i)] $M_0\cong N$;
\item[(ii)] $M_n=M$;
\item[(iii)] for $k\in\{1,2,\dots,m\}$, $r(M_k)-r(M_{k-1})=1$ and $|E(M_k)-E(M_{k-1})|\le 3$; and
\item[(iv)] for $m<k\le n$, $r(M_k)=r(M)$ and $|E(M_k)-E(M_{k-1})|=1$.
\end{enumerate}
Moreover, when $|E(M_k)-E(M_{k-1})|=3$, for some $1\le k\le m$, $E(M_k)-E(M_{k-1})$ is a triad of $M_k$.}
\bigskip

The Strong Splitter Theorem says that we can obtain up to isomorphism $M$ starting with $N$ and at each step doing a 3-connected single-element extension or coextension, such that at most two consecutive single-element extensions may occur in the sequence before a single-element coextension must occur, unless the rank of the minors involved are the same as the rank of $M$. Moreover, as the last line indicates, if two consecutive single-element extensions by elements $\{e_1, e_2\}$ are followed by a coextension by element $f$, then $\{e_1, e_2, f\}$ forms a triad in the resulting matroid. 

Motivated by this result, observe that every rank-$r$ extremal matroid has a  rank-$r$  root matroid $\alpha _r$ that is a cosimple single-element coextension of the rank-$(r-1)$ root, $\alpha _{r-1}$ or of its simple single-element extensions, or of its simple double-element extensions (formed in the one specific manner described by the Strong Splitter Theorem.) There is no reason to asssume {\it a priori} that the rank-$r$ root matroid is unique. The rank-$r$ root matroid(s) could have 1, 2, or 3 more elements than the rank-$(r-1)$ root(s), and in each case there could be many non-isomorphic rank-$r$ roots.  However, when $\alpha _r$ happens to be  unique, we get a recursive way of defining it and then in turn a recursive way of defining the rank-$r$ extremal matroid $\Omega _r$. In subsequent papers we will develop strategies for ensuring uniqueness of the root matroid in certain types of excluded minor classes.

For $EX[P_9^*]$ the rank-$r$ root matroid   $\alpha_r$ happens to be unique and it has $(3r-5)$ elements. A matrix representation is shown below: 

\begin{figure}[h]
\centering
\epsfxsize 6in \epsfbox{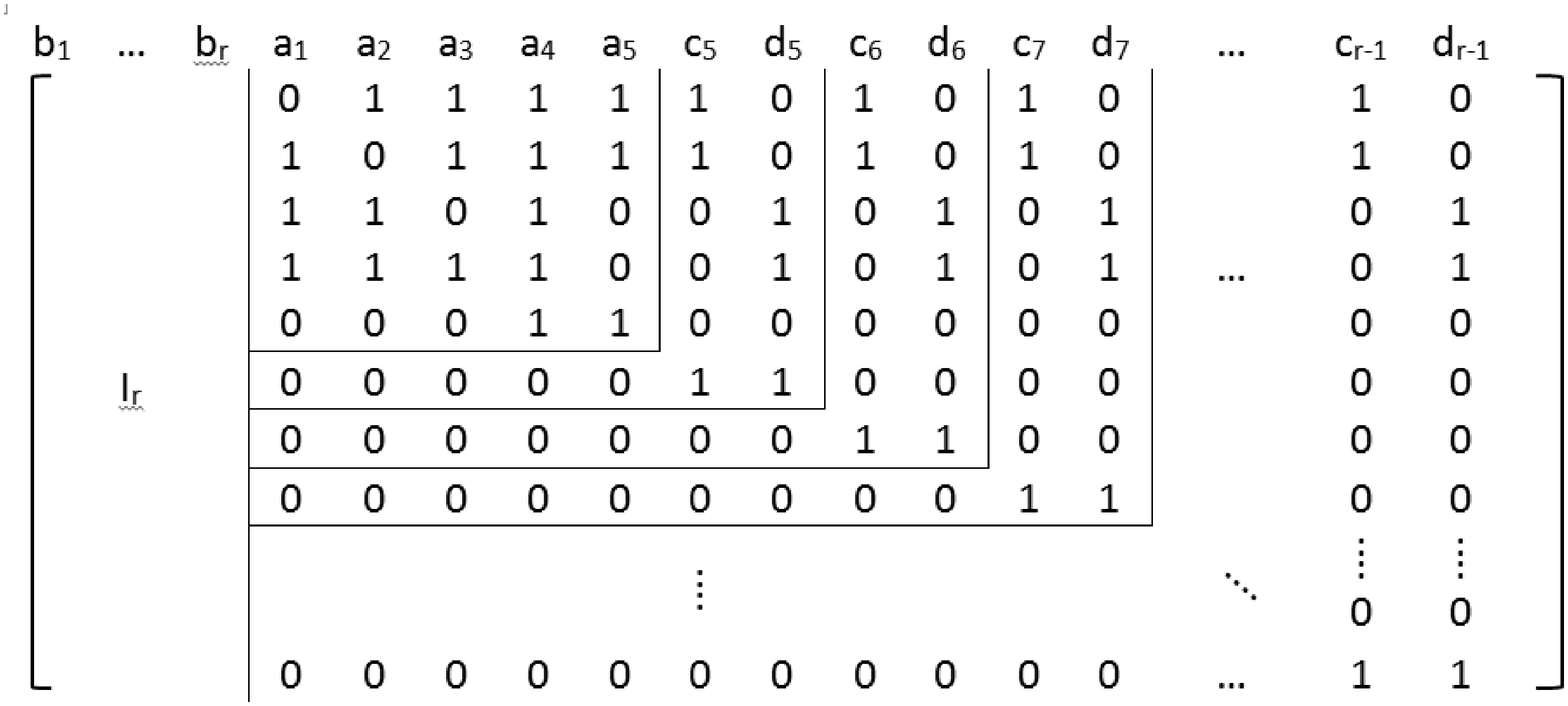}
\caption{The rank-$r$ $(3r-5)$-element root matroid $\alpha _r$, for $r\ge 5$, in $EX[P_9^*]$}
\end{figure}

The above matrix representation reveals how $\alpha _r$ is recursively constructed from $\alpha _5$ shown below, which is the starting matroid for this family:
\small
\[ 
\alpha _5=\left[ 
\begin{array}{c|ccccc}
&    0&1&1&1&1 \\
&    1&0&1&1&1\\
I_5& 1&1&0&1&0 \\
&    1&1&1&1&0 \\
&    0&0&0&1&1
\end{array} 
\right] 
\]
\normalsize
\noindent 
Observe that $\alpha _6$ (shown below) is obtained from $\alpha _5$ by adding two columns $c_5=[1 1 0 0 0 ]^T$ and $d_5=[0 0 1 1 0 ]^T$ and lifting by row $[0 0 0 0 0 1 1]$. In other words add $c_5$ and $d_5$ to form triangles with $\{b_1, b_2, c_5\}$ and $\{b_3, b_4, d_5\}$, respectively; then lift elements $c_5$ and $d_5$ into the next dimension to form a triad $\{c_5, d_5, b_6\}$ with the new lift element $b_6$. 
\small
\[
\alpha _6=\left[ 
\begin{array}{c|ccccccc}
&    0&1&1&1&1&1&0 \\
&    1&0&1&1&1&1&0\\
I_6& 1&1&0&1&0&0&1 \\
&    1&1&1&1&0&0&1 \\
&    0&0&0&1&1&0&0 \\
&    0&0&0&0&0&1&1
\end{array} 
\right] 
\] 
\normalsize

In general, $\alpha _r$ is formed as follows: add parallel elements $\{c_5, c_6,  \dots , c_{r-1}\}$ so that each forms a triangle with basis points $b_1$ and $b_2$; add parallel elements $\{d_5, d_6,  \dots , d_{r-1}\}$ so that each forms a triangle with basis points $b_3$ and $b_5$; do a sequence of lifts by adding new basis elements $\{b_6, \dots , b_r\}$ to form triads $\{b_i, c_{i-1}, d_{i-1}\}$ for $i=6, \dots , r$. 

Once the construction of the rank-$r$ root matroid is understood, the construction of the rank $r$-extremal matroid folows easily. To obtain $\Omega _5$ from $\alpha _5$ add five columns $c_5=[1 1 0 0 0]^T$, $d_6=[0 0 1 1 0]^T$, $e_6=[1 1 1 0 0]^T$, $f_6=[0 0 1 1 1]^T$, $g_{6, 1}=[1 1 1 1 0 0]^T$.
To obtain $\Omega _6$ from $\alpha _5$ add six columns $c_6=[1 1 0 0 0 0]$, $d_6=[0 0 1 1 0 0]^T$, $e_6=[1 1 1 0 0 0]^T$, $f_6=[0 0 1 1 1 0]^T$, $g_{6, 1}=[1 1 1 1 0 0]^T$ and $g_{6, 2}=[1 1 1 1 0 1]^T$. Matrix representations for $\Omega _5$ and $\Omega _6$ are shown below:

\small
\[ 
\Omega _5=\left[ 
\begin{array}{c|cccccccccc}
&    0&1&1&1&1&1&0&1&0&1 \\
&    1&0&1&1&1&1&0&1&0&1\\
I_5& 1&1&0&1&0&0&1&1&1&1 \\
&    1&1&1&1&0&0&1&0&1&1 \\
&    0&0&0&1&1&0&0&0&1&0
\end{array} 
\right] 
\Omega _6=\left[ 
\begin{array}{c|ccccccccccccc}
&    0&1&1&1&1&1&0&1&0&1&0&1&1 \\
&    1&0&1&1&1&1&0&1&0&1&0&1&1\\
I_6& 1&1&0&1&0&0&1&0&1&1&1&1&1 \\
&    1&1&1&1&0&0&1&0&1&0&1&1&1 \\
&    0&0&0&1&1&0&0&0&0&0&1&0&0 \\
&    0&0&0&0&0&1&1&0&0&0&0&0&1
\end{array} 
\right] 
\] 
\normalsize

\noindent In general to obtain $\Omega _r$ from $\alpha _r$ add $r$ columns $c_r=[1 1 0 0 0 0,  \dots , 0 0]^T$, $d_r=[0 0 1 1 0 0,  \dots , 0 0]^T$, $e_r=[1 1 1 0 0 0, \dots , 0 0]^T$, $f_r=[0 0 1 1 1 0  \dots 0 0]^T$, and
$g_{r, 1}=[1 1 1 1 0 0 0 0, \dots , 0 0 0]^T$, 
$g_{2, 2}=[1 1 1 1 0 1 0 0, \dots , 0 0 0]^T$, 
$g_{2, 3}=[1 1 1 1 0 0 1 0, \dots , 0 0 0]^T$up to 
$g_{2, r-5}=[1 1 1 1 0 0 0, \dots , 0 1 0]^T$
$g_{r, r-4}=[1 1 1 1 0 0 0, \dots , 0 0 1]^T$. A matrix representation for $\Omega _r$ is shown below.

\begin{figure}[h]
\centering
\epsfxsize 7in \epsfbox{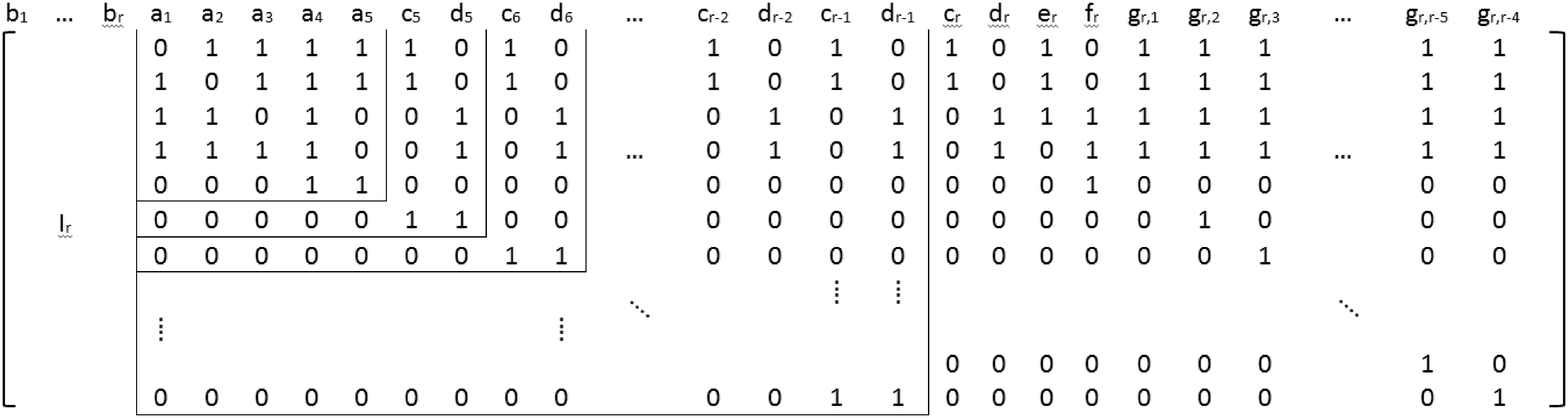}
\caption{The rank-$r$ $(4r-5)$-element extremal matroid $\Omega _r$, for $r\ge 5$, in $EX[P_9^*]$}
\end{figure}

\bigskip
\noindent {\bf Proposition 2.1.} {\it The matroid $\Omega _r$ has no $P_9^*$-minor.}
\bigskip

\noindent {\bf Proof.} Observe that, the matroid $\alpha _5$ has no odd element circuits, whereas $P_9^*$ has odd element circuits. Thus the 10-element matroid $\alpha _5$ cannot have the 9-element matroid $P_9^*$ as a minor (also shown in Appendix Table 1b). Since $\alpha _r$ is constructed by adding triangles and triads in the recursive manner  described above, it also has no odd circuits. (Deleting any of the columns $c_5, d_5, \dots c_r, d_r$ gives a matroid that is not 3-connected.) Next, using the facts that $\Omega _5$ has no $P_9^*$-minor (Table 2a in the Appendix) and $\Omega _r/b_r\backslash c_r, d_r, g_{r, r-4} = \Omega _{r-1}$, an induction argument shows that $\Omega _r$ has no $P_9^*$-minor. $\qed$

\bigskip

\section {Extremal matroids in $EX[P_9^*]$}

The main result in this section is given below.
\bigskip

\noindent{\bf Theorem 3.1}
\footnote {Note to referee: A similar result was independently proven by Ding and Wu. The entire story is attached for the editor. All email correspondence on this paper is included with their permission largely to establish that there is no disagreement of any sort and we have plans to work together in future. In this instance we proved the same result using completely different methods. It is up to the editor if he thinks the referee needs to see the emails. Ding and Wu use completely different techniques and their result is not a growth rate result. Moreover, their result has a computer component.}
{\it A binary $3$-connected  non-regular matroid $M$ has no $P_9^*$-minor if and only if $M$ is isomorphic to $F_7$, $PG(3,2)$, $R_{16}$, $Z_r$ for $r\ge 4$, $\Omega_r$ for $r\ge 5$, or one of their $3$-connected deletion-minors.  }
\bigskip

 Since $P_9$ is a rank-5 matroid, the rank-3 binary extremal matroid $F_7$ and the rank-4 binary extremal matroid $PG(3,2)$ are automatically in $EX[P_9^*]$. The class also contains a rank-5 16-element extremal matroid $R_{16}$. These three finite matroids are internally 4-connected. Matrix representations for them  are given below:
\small
\[
F_7=\left[ 
\begin{array}{c|cccc}
&   0&1&1&1 \\
I_3&1&0&1&1 \\
&   1&1&0&1 \\
\end{array} 
\right] 
PG(3, 2)=\left[ 
\begin{array}{c|cccccccccccc}
&   0&0&0&0&1&1&1&1&1&1&1  \\
I_4& 0&1&1&1&0&0&0&1&1&1&1 \\
&   1&0&1&1&0&1&1&0&0&1&1 \\
&   1&1&0&1&1&0&1&0&1&0&1
\end{array} 
\right] 
\]
\normalsize
\small
\[
R_{16}=\left[ 
\begin{array}{c|ccccccccccc}
&    1&0&0&1&1&0&0&1&1&1&1 \\
&    1&1&0&0&1&1&1&0&0&1&1 \\
I_5& 1&1&1&0&0&0&1&1&1&0&1 \\
&    0&1&1&1&0&1&0&0&1&1&1 \\
&    0&0&1&1&1&1&1&1&0&0&1
\end{array}
\right] 
\] 
\normalsize

In [\ref{Oxley1987}] Oxley proved that a $3$-connected binary non-regular matroid $M$ has no minor isomorphic to $P_9$ or $P_9^*$ if and only if $M$ is isomorphic to $F_7$, $F_7^*$, $Z_r$, $Z_r^*$, $Z_r\backslash b_r$, or  $Z_r\backslash c_r$,  for some $r\ge 4$. The matroid $Z_r$ can be represented by the matrix $[I_r|D]$, where $D$ has $r+1$ columns with zeros on the diagonal and ones elsewhere. Label the columns $b_1, \dots , b_r, a_1, \dots , a_r, c_r$. The matroids $Z_r\backslash b_r$ and $Z_r\backslash c_r$ are 3-connected single-element deletion-minors of $Z_r$ and $Z_r\backslash \{b_r, c_r\}= Z_{r-1}^*$. Moreover, $Z_r\backslash b_r$ and $Z_r\backslash c_r$  are self-dual. Thus, the infinite family of extremal matroids for $EX[P_9, P_9^*]$ is $Z_r$, for $r\ge 4$. Observe that, the rank-$r$ root matroid of $Z_r$ is $Z_{r-1}^*$.

\bigskip
\noindent {\bf Proof.} Suppose $M$ is a 3-connected non-regular matroid in $EX[P_9^*]$. Then $M$ has a minor isomorphic to $F_7$ or $F_7^*$ [\ref{Oxley2012}, 10.1.2].  Observe that $F_7=PG(2,2)$ and as such has no extensions in the class of binary matroids. Coextensions of $F_7$ are duals of extensions of $F_7^*$. Thus we may focus on the extensions of $F_7^*$.  By the Splitter Theorem, $M$ has a sequence of 3-connected minors $M_0\cong F_7^*$, $M_1$, $M_2$, \dots $M_k=M$ such that each matroid in the sequence is a single-element extension or coextension of the previous one. By the Strong Splitter Theorem we may do at most two extensions before a coextension occurs. Moreover, if a matroid is 3-connected, then its simple single-element extensions and cosimple single-element coextensions are also 3-connected. The next lemma is the first of three lemmas in the proof.

\bigskip

\noindent {\bf Lemma 2.1.}
\footnote{Note to referee: This lemma is the main theorem in [\ref{Oxley1987}] re-proved here to showcase a simple application of the Strong Splitter Theorem, as a preview to the more complicated application in this paper. The proof can be removed.}
 Suppose $M$ is a $3$-connected binary non-regular matroid with no $P_9$ nor $P_9^*$-minor. Then the rank-$r$ root matroid is $Z_{r-1}^*$ and the extremal matroid is $Z_r$.
\bigskip

\noindent {\bf Proof.} The proof is by induction on $r\ge 4$. Suppose $r=4$. Observe that $AG(3, 2)$ and $S_{8}$ are the two non-isomorphic simple single-element extensions of $F_7^*$. Since they are self-dual, they are also the coextensions of $F_7$. The matroid $S_8$ has two non-isomorphic simple single-element extensions $P_9$ and $Z_4$ and $AG(3, 2)$ has one simple single-element extension $Z_4$.  Observe that $Z_4$ is the rank-4 extremal matroid in $EX[P_9, 
P_9^*]$. Moreover, $AG(3,2)=Z_4\backslash c_4$ and $S_8=Z_4\backslash a_4$ and both are self-dual. Their coextension is $Z_4^*$, which becomes the rank-5 root matroid. Thus the result holds for $r=4$.
Assume the result holds for all non-regular matroids of rank at most $(r-1)$  in $EX[P_9, P_9^*]$ for $r\ge 5$. Suppose $M$ is a rank-$r$ matroid in $EX[P_9, P_9^*]$. 
\bigskip

\noindent {\bf Claim A.} The rank-$(r-1)$ root $Z_{r-2}^*$  gives rise to the rank-$r$ root $Z_{r-1}^*$.
\bigskip
\bigskip

\noindent {\bf Proof.} By the Strong Splitter Theorem, $M$ must be a cosimple single-element coextension of  $Z_{r-2}^*$, $Z_{r-1}\backslash c_{r-1}$, $Z_{r-1}\backslash a_{r-1}$, or $Z_{r-1}$. Moreover, if $M$ is a cosimple single-element coextension of $Z_{r-1}$, then $\{b_r, a_r, c_r\}$ forms a triad in $M$.

By the induction hypothesis the only rows that can be added to $Z_{r-3}$   are $[1 1 \dots 1 0]$ and $[1 1 \dots 1 1]$. Adding $[1 1 \dots 1 0]$ gives $Z_{r-2}\backslash c_{r-2}$ and adding $[1 1 \dots 1 1]$ gives  $Z_{r-2}\backslash a_{r-2}$. Adding both gives $Z_{r-2}^*$. Therefore $Z_{r-2}^*$ has no further cosimple coextensions in $EX[P_9, P_9^*]$. 

The only simple single-element extensions of $Z_{r-2}^*$ in $EX[P_9, P_9^*]$ are obtained by adding columns $a_{r-1}=[1 1 \dots 1 0]^T$ and $c_{r-1}=[1 1 \dots 1 1]^T$ giving respectively, $Z_{r-1}\backslash c_{r-1}$ and $Z_{r-1}\backslash a_{r-1}$. However, $Z_{r-1}\backslash c_{r-1}$ and $Z_{r-1}\backslash a_{r-1}$  are also single-element coextensions of $Z_{r-2}$ by rows $[1 1  \dots   1 0]$ and $[1 1  \dots   1 1]$, respectively. Adding both these rows to $Z_{r-2}$ gives $Z_{r-1}^*$. 

Adding to $Z_{r-2}^*$ both columns $c_{r-2}$ and $a_{r-2}$ gives $Z_{r-1}$. Lastly, the only cosimple single-element coextension of $Z_{r-1}$ we must check is the matroid $Z_{r-1}'$ formed by adding row $[0 0 \dots 0 1 1]$.  Observe that, $Z_{r-1}'/\{b_5, b_6, \dots b_{r-1}\}\backslash \{a_5, a_6, \dots a_{r-1}\}$ is the matroid $M$ shown below and $M\backslash 6\cong P_9^*$.
\small 
\[ 
M=\left[ 
\begin{array}{c|ccccc}
&    0&1&1&1&1  \\
&    1&0&1&1&1 \\
I_5& 1&1&0&1&1  \\
&    1&1&1&0&1  \\
&    0&0&0&1&1 
\end{array} 
\right] 
\]
\normalsize

\noindent Therefore, the rank-$r$  root matroid is $Z_{r-1}^*$.
 
\bigskip

\noindent {\bf Claim B.} The rank-$r$ root $Z_{r-1}^*$  extends to the rank-$r$ extremal matroid $Z_r$.
\bigskip

\noindent {\bf Proof.}  We will prove that the only columns that can be added to $Z_{r-1}^*$ are $c_r=[1 1 \dots 1 1]^T$ and $a_r=[1 1 \dots 1 0]^T$. 
To prove this, first observe that $Z_{r-1}^*/b_r=Z_{r-1}\backslash c_{r-1}$ and $Z_{r-1}^*/b_{r-1}=Z_{r-1}\backslash a_{r-1}$. By the induction hypothesis applied to $Z_{r-1}^*/b_r$ the only columns that can be added  are $c_{r-1}$ with a zero or one in the last position, $b_1, b_2, \dots b_{r-2}, b_{r-1}$ with a one in the last position, and $a_1, a_2, \dots a_{r-2}, a_{r-1}$ with the entry in the last position switched. They are:
\medskip

 \begin{enumerate}
\item {\bf Type I columns}: $c_{r-1}^0 = [1 1  \dots 1 0]^T$ and $c_{r-1}^1 = [1 1  \dots 1 1]^T$;
\item {\bf Type II columns}: $b_1^1=[1 0 0 \dots 0 1]^T$, $b_2^1=[0 1 0 \dots 0 1]^T$ up to $b_{r-2}^1=[ 0 0 0 \dots 0 1 0 1]^T$, $b_{r-1}^1=[0 0 0 \dots 0 1 1]^T$; and 
\item {\bf Type III columns}: $a_1^0 = [0 1 1 1 \dots 1 1 1 0]^T$, $a_2^0 = [1 0 1 1 \dots 1 1 1 0]^T$, up to $a_{r-2}^0 = [1 1 1 \dots 1 0 1 0]^T$, $a_{r-1}^0 = [1 1 1 \dots 1 1 0 0]^T$.
\end{enumerate}

\medskip
\noindent Similary, the only columns that can be added to $Z_{r-1}^*/b_{r-1}$ are $a_{r-1}$ with a zero or one, $b_1, b_2, \dots b_{r-2}, b_r$ with a one in the second-last position, and $a_1, a_2, \dots a_{r-2}, c_{r-1}$ with the entry in the second-last position switched. They are:
\medskip

 \begin{enumerate}
\item {\bf Type I columns}: $a_{r-1}^0 = [1 1  \dots 0 0]^T$ and $a_{r-1}^1 = [1 1  \dots 1 0]^T$;
\item {\bf Type II columns}: $b_1^1=[1 0 0 \dots 1 0]^T$, $b_2^1=[0 1 0 \dots 1 0]^T$ to $b_{r-2}^1=[ 0 0 0 \dots 0 1 1 0]^T$, $b_r^1=[0 0 0 \dots 0 1 1]^T$; and 
\item {\bf Type III columns}: $a_1^0 = [0 1 1 1 \dots 1 1 0 1]^T$, $a_2^0 = [1 0 1 1 \dots 1 1 0 1]^T$, to $a_{r-2}^0 = [1 1 1 \dots 1 0 0 1]^T$, and $a_{r-1}^1 = [1 1 1 \dots 1 1 1 1]^T$.
\end{enumerate}

\medskip

\noindent Observe that the only overlapping columns are $[1 1  \dots 1 0]^T$, $[1 1  \dots 1 1]^T$, and $[0 0 \dots 0 1 1]$. The first is $a_r$ and the second is $c_r$. They give the single-element extensions $Z_r\backslash c_r$ and $Z_r\backslash a_r$ and together the double-element extension $Z_r$. 

Lastly, let $Z_{r-1}^*+b_r^1$ be the matroid obtained by adding $b_r^1=[0 0 \dots 1 1]$ to $Z_{r-1}^*$. Observe that $$(Z_{r-1}^*+b_{r-1}^1) / \{b_4, \dots , b_{r-2}\} \backslash \{a_4, \dots , a_{r-2}\} = Z_5^*+b_4^1.$$ The matroid $Z_5^*+b_4^1$ shown below has a $P_9^*$-minor.

\small
\[ 
Z_5^*+b_4^1=\left[ 
\begin{array}{c|ccccc}
&    0&1&1&1&0 \\
&    1&0&1&1&0 \\
I_5& 1&1&0&1&0 \\
&    1&1&1&0&1 \\
&    1&1&1&1&1
\end{array} 
\right] 
\] 
\normalsize 

\noindent Thus, $Z_r$ is the rank-$r$ extremal matroid. 
\bigskip

This completes the proof of Lemma 2.1. $\qed$

\bigskip
 
Returning to the proof of Theorem 1.1, suppose $M$ is a $3$-connected binary  non-regular matroid with a $P_9$-minor, but no $P_9^*$-minor. 
From Tables 1a and 1b in the Appendix we see that $P_9$ has three non-isomorphic simple single-element extensions, $D_1$, $D_2$, and $D_3$, and eight non-isomorphic cosimple single-element coextensions of which just one matroid $E_7$ has no $P_9^*$-minor. 
 
\small
\[ 
E_7=\alpha_5=\left[ 
\begin{array}{c|ccccc}
&    0&1&1&1&1 \\
&    1&0&1&1&1 \\
I_5& 1&1&0&1&0 \\
&    1&1&1&1&0 \\
&    0&0&0&1&1
\end{array} 
\right] 
\] 
\normalsize

\noindent Since $P_9^*$ is a rank 5 matroid, $EX[P_9^*]$ contains $PG(3,2)$. We will prove that a matroid of rank $r\ge 5$ in $EX[P_9^*]$ has an $E_7$-minor. Note that $E_7=\alpha_5$.

\bigskip
\noindent {\bf Lemma 2.2.} 
\footnote{Note to referee: Lemma 2.2 and 2.3 follow from the main result in [\ref{KinganSubmitted}], however, the details of that more general result are not necessary as stand-alone proofs of Lemmas 2.2 and 2.3 are very short comparitively and including it makes the paper self-contained and easier to read for the referee. The proofs of Lemmas 2.2 and 2.3 and Table 4 can be removed if the referee thinks so. }
Suppose $M$ is a 3-connected binary matroid with a $P_9$-minor and no $P_9^*$ minor and rank at least 5. Then $M$ has an $E_7$-minor.
\bigskip

\noindent {\bf Proof.} By the Strong Splitter Theorem $M$ must be a cosimple single-element coextension of $P_9$ or of its single-element extensions $D_1$, $D_2$, or $D_3$, or of its double-element extensions $X_1$, $X_2$, $X_3$ formed with row $[0 0 0 0 0 1 1]$. Observe that if $M$ is a coextension of $P_9$, then as  mentioned earlier $M\cong E_7$ and $E_7$ is formed by adding just one row $[0 0 0 1 1]$. 

Suppose $M$ is a coextension of $D_1$, $D_2$, or $D_3$ shown below:

\tiny
\[ 
D_1=\left[ 
\begin{array}{c|cccccc}
&   0&1&1&1&1&1 \\
I_4&1&0&1&1&1&1 \\
&   1&1&0&1&0&1 \\
&   1&1&1&1&0&0
\end{array} 
\right] 
D_2=\left[ 
\begin{array}{c|cccccc}
&   0&1&1&1&1&1 \\
I_4&1&0&1&1&1&0 \\
&   1&1&0&1&0&0 \\
&   1&1&1&1&0&1
\end{array} 
\right] 
D_3=\left[ 
\begin{array}{c|cccccc}
&   0&1&1&1&1&0 \\
I_4&1&0&1&1&1&0 \\
&   1&1&0&1&0&1 \\
&   1&1&1&1&0&1
\end{array} 
\right] 
\] 
\normalsize

\noindent Type I rows that can be added to $D_1$, $D_2$, and $D_3$ are $[0 0 0 1 1 0]$ and $[0 0 0 1 1 1]$. Type II rows are $[1 0 0 0 0 1]$, $[0 1 0 0 0 1]$, $[0 0 1 0 0 1]$, and $[0 0 0 0 1 1]$. Type III rows are the rows of $D_1$, $D_2$, and $D_3$ with the last position switched. Table 4 in the Appendix shows that in all cases the resulting coextension has a $P_9^*$-minor or an $E_7$-minor.
Lastly, the matrices $X_1$, $X_2$, and $X_3$ with row $[0 0 0 0 0 1 1]$ are shown below. They have a $P_9^*$-minor or an $E_7$-minor. $\qed$

 \tiny
\[ 
X_1'=\left[ 
\begin{array}{c|ccccccc}
&   0&1&1&1&1&1&1\\
I_5&1&0&1&1&1&0&0 \\
&   1&1&0&1&0&0&1 \\
&   1&1&1&1&0&1&0 \\
&   0&0&0&0&0&1&1
\end{array} 
\right] 
X_2'=\left[ 
\begin{array}{c|ccccccc}
&   0&1&1&1&1&1&0\\
I_5&1&0&1&1&1&1&0 \\
&   1&1&0&1&0&1&1 \\
&   1&1&1&1&0&0&1 \\
&   0&0&0&0&0&1&1
\end{array} 
\right] 
X_3'=\left[ 
\begin{array}{c|ccccccc}
&   0&1&1&1&1&1&0\\
I_5&1&0&1&1&1&0&0 \\
&   1&1&0&1&0&0&1 \\
&   1&1&1&1&0&1&1 \\
&   0&0&0&0&0&1&1
\end{array} 
\right] 
\] 
\normalsize

\bigskip

The simple single-element extensions and cosimple single-element coextensions of $E_7$ are shown in Tables 2a and 2b in the Appendix, respectively. All three single-element coextensions of $E_7$ have a $P_9^*$-minor, but there are three single-element extensions with no $P_9^*$-minor (extensions 2, 3, and 5). Let $\alpha_{5,1}=(\alpha_5, ext2)$, $\alpha_{5,2}=(\alpha_5, ext3)$ and  $\alpha_{5,3}=(\alpha_5, ext5)$. Matrix representations for $\alpha_{5,1}$, $\alpha_{5,2}$, and $\alpha_{5,3}$ are shown below.

\tiny
\[
\alpha_{5,1}=\left[ 
\begin{array}{c|cccccc}
&    0&1&1&1&1&1 \\
&    1&0&1&1&1&1\\
I_5& 1&1&0&1&0&0 \\
&    1&1&1&1&0&0 \\
&    0&0&1&1&1&0
\end{array} 
\right]
\alpha_{5,2}=\left[ 
\begin{array}{c|cccccc}
&    0&1&1&1&1&1 \\
&    1&0&1&1&1&1\\
I_5& 1&1&0&1&0&1 \\
&    1&1&1&1&0&0 \\
&    0&0&1&1&1&0
\end{array} 
\right]
\alpha_{5, 3}=\left[ 
\begin{array}{c|cccccc}
&    0&1&1&1&1&0 \\
&    1&0&1&1&1&1\\
I_5& 1&1&0&1&0&0 \\
&    1&1&1&1&0&1 \\
&    0&0&1&1&1&1
\end{array} 
\right] 
\]

\normalsize

\bigskip
\noindent {\bf Lemma 2.3.} {\it Suppose $M$ is a $3$-connected binary non-regular matroid with an $\alpha_{5, 3}$-minor and no  $P_9^*$-minor. Then $r(M)\le 5$.   }
\bigskip

\noindent {\bf Proof.} Observe from Table 2a that $\alpha_{5, 1}$ is obtained by adding columns $a=[ 0 0 1 1 0 ]$, 
$b=[1 1 0 0 0]$,   and $c=[1 1 1 1 0]$; $\alpha_{5, 2}$ is obtained by adding column  $d=[0 0 1 1 1]$ and $e=[1 1 1 0 0]$; and  $\alpha_{5, 3}$ is obtained by adding column  $f=[0 1 0 1 1]$, $g=[0 1 1 0 1]$, $h=[1 0 0 1 1]$, and $i=[1 0 1 0 1]$. We can check that $\alpha_{5, 3}$ (formed by adding column $f$ to $E_7$) has only one simple single-element extension in $EX[P_9^*]$ and it is obtained by adding column $g$, $h$, $i$, $d$, and $e$. Up to isomorphism all five columns give the same single-element extension. Let us call this matroid $\alpha_{5, 3, 1}$ obtained by adding say column $g$.  

\small
\[
\alpha_{5, 3, 1}=\left[ 
\begin{array}{c|ccccccc}
&    0&1&1&1&1&0&0 \\
&    1&0&1&1&1&1&1\\
I_5& 1&1&0&1&0&0&1 \\
&    1&1&1&1&0&1&0 \\
&    0&0&0&1&1&1&1
\end{array} 
\right] 
\]  
\normalsize

Similarly, adding to $\alpha_{5, 3, 1}$ any one of columns $h, i, d, e$, say $h$, gives $\alpha_{5, 3, 1, 1}$, and so on; we get $\alpha_{5, 3, 1, 1}$ by adding $i$; $\alpha_{5, 3, 1, 1, 1}$ by adding $d$; and finally $\alpha_{5, 3, 1, 1, 1, 1}=R_{16}$ by adding $e$. Thus, we  determined all the rank 5 members in $EX[P_9^*]$ with an $\alpha_{5, 3}$-minor. 

It remains to show that $r(M)\le 5$. Suppose $M$ is a rank 6 cosimple single-element coextension of a rank 5 matroid in $EX[P_9^*]$. By the Strong Splitter Theorem, $M$ is a cosimple single-element coextensions of $E_7=\alpha _5$ or $\alpha_{5, 3}$ or $\alpha _{5, 3, 1}$ (with row $[0 0 0 0 0 1 1]$). Observe that $E_7$ has no cosimple coextensions in $EX[P_9^*]$. Therefore, $M$ may be a cosimple single-element coextension of $\alpha_{5, 3}$. We may add to $\alpha_{5, 3}$ only Type II and Type III rows. 
Type II rows are $[1 0 0 0 0 1]$, $[0 1 0 0 0 1]$, $[0 0 1 0 0 1]$, $[0 0 0 1 0 1]$, and $[0 0 0 0 1 1]$. Type III rows are $[0 1 1 1 1 1]$, $[1 0 1 1 1 0]$, $[1 1 0 1 0 1]$, $[1 1 1 1 0 0]$, and $[0 0 1 1 1 0]$.  The cosimple single-element coextensions of $\alpha_{5, 3}$ are shown in Table 3 with these rows highlighted. Observe that, all the cosimple single-element coextensions of $\alpha_{5, 3}$ have a $P_9^*$-minor. 
Lastly, $\alpha _{5, 3, 1}$ with row $[0 0 0 0 0 1 1]$ is shown below:

\small
\[
\alpha_{5, 3, 1}'=\left[ 
\begin{array}{c|ccccccc}
&    0&1&1&1&1&0&0 \\
&    1&0&1&1&1&1&1 \\
I_6& 1&1&0&1&0&0&1 \\
&    1&1&1&1&0&1&0 \\
&    0&0&0&1&1&1&1 \\
&    0&0&0&0&0&1&1
\end{array} 
\right] 
\]  
\normalsize
\noindent Observe that $\alpha_{5, 3, 1}'/3\backslash \{7, 8, 10\}\cong P_9^*$. Thus, we conclude that $r(M)\le 5$. $\qed$

\bigskip
\noindent {\bf Lemma 2.4.} \footnote{Note to referee: As mentioned earlier, the proofs of Lemmas 2.1, 2.2, and 2.3 may be removed and replaced with a one page explanation citing previous results. Lemma 2.4 is a new result.} {\it Suppose $M$ is a $3$-connected binary matroid with an $\alpha_{5, 1}$- or $\alpha_{5, 2}$-minor and no  $P_9^*$-minor. Then the rank-$r$ root matroid is $\alpha_r$ and the extremal matroid is $\Omega_r$. }
\bigskip

\noindent {\bf Proof.} The proof is by induction on $r\ge 5$. Suppose $r=5$. As noted earlier, Table 2a in the Appendix shows that $\alpha _{5,1}$ is formed by adding columns $c_5=[1 1 0 0 0]^T$, $d_5=[0 0 1 1 0]^T$, or $g_{5, 1}=[1 1 1 1 0]^T$ and $\alpha_{5, 2}$ is formed by adding columns $e_5=[1 1 1 0 0]^T$ or $f_5=[0 0 1 1 1]^T$. Adding all these columns to $\alpha _5$ gives $\Omega _5$. Table 2b shows that every cosimple single-element coextension of $\alpha_5$ has a $P_9^*$-minor and (using the same method as for $\alpha _{5,3}$) Table 3 shows that every cosimple single-element coextension of $\alpha _{5, 1}$ and $\alpha _{5, 2}$ also have a $P_9^*$-minor. 

It is easy to check that $\alpha _{5, 1}$ (with $c_5$) has two simple single-element extensions in $EX[P_9^*]$, namely, $\alpha _{5, 1, 1}$, formed by adding $d_5$ or $g_{5, 1}$, and $\alpha _{5, 1, 2}$ formed by adding  $e_5$ or $f_5$. The matroid $\alpha _{5, 2}$ (with $e_5$) also has two single-element extensions, $\alpha _{5, 2, 1}$ formed by adding $c_5$, $d_5$ or $g_{5,1}$ and $\alpha _{5, 2, 2}$ formed by adding $f_5$. Further, note that  $\alpha _{5, 2, 1}=\alpha _{5, 1, 2}$. 

\tiny
\[ 
\alpha_{5, 1, 1}=\left[ 
\begin{array}{c|ccccccc}
&    0&1&1&1&1&1&0 \\
&    1&0&1&1&1&1&0 \\
I_5& 1&1&0&1&0&0&1 \\
&    1&1&1&1&0&0&1 \\
&    0&0&0&1&1&0&0
\end{array} 
\right] 
\alpha_{5, 1, 2}=\left[ 
\begin{array}{c|ccccccc}
&    0&1&1&1&1&1&1 \\
&    1&0&1&1&1&1&1 \\
I_5& 1&1&0&1&0&0&1 \\
&    1&1&1&1&0&0&0 \\
&    0&0&0&1&1&0&0
\end{array} 
\right] 
\alpha_{5, 2, 2}=\left[ 
\begin{array}{c|ccccccc}
&    0&1&1&1&1&0 \\
&    1&0&1&1&1&0 \\
I_5& 1&1&0&1&1&1 \\
&    1&1&1&1&0&1 \\
&    1&1&0&0&0&1
\end{array} 
\right] 
\] 
\normalsize

By the Strong Splitter Theorem we must only check one single-element coextension of  $\alpha _{5, 1, 1}$, $\alpha _{5, 1, 2}$, and $\alpha _{5, 2, 2}$, namely the one formed by adding row $[0 0 0 0 0 1 1]$. Observe that $\alpha _{5, 1, 1}$ with row $[0 0 0 0 0 1 1]$ is precisely $\alpha _6$.  The matroid $\alpha _{5, 1, 2}'$ and $\alpha _{5, 2, 2}'$ obtained by adding row $[0 0 0 0 0 1 1]$ to $\alpha _{5, 1, 2}$ and $\alpha _{5, 2, 2}$, respective are shown below:

\small
\[
\alpha_{5, 1, 2}'=\left[ 
\begin{array}{c|ccccccc}
&    0&1&1&1&1&1&1 \\
&    1&0&1&1&1&1&1 \\
I_6& 1&1&0&1&0&0&1 \\
&    1&1&1&1&0&0&0 \\
&    0&0&0&1&1&0&0 \\
&    0&0&0&0&0&1&1 
\end{array}
\right]  
\alpha_{5, 2, 2}'=\left[ 
\begin{array}{c|ccccccc}
&    0&1&1&1&1&0 \\
&    1&0&1&1&1&0 \\
I_5& 1&1&0&1&1&1 \\
&    1&1&1&1&0&1 \\
&    1&1&0&0&0&1 \\
&    0&0&0&0&1&1 
\end{array} 
\right] 
\] 
\normalsize

\noindent Observe that $\alpha_{5, 1, 2}' / 5 \backslash \{9, 10, 11\}\cong P_9^*$ and $\alpha_{5, 2, 2}' / 5 \backslash \{8, 10, 11\}\cong P_9^*$. Thus, $M = \alpha _{5, 1, 2}' = \alpha _6$ and the result holds for $r=5$.

Assume the result holds for all non-regular matroids of rank at most $(r-1)$ in $EX[P_9^*]$ with a minor isomorphic to $\alpha_{5,1}$ or $\alpha_{5,2}$. Suppose $M$ is a rank-$r$  matroid in $EX[P_9^*]$ with a minor isomorphic to $\alpha_{5,1}$ or $\alpha_{5,2}$. We will first show that $\alpha _{r-1}$ gives rise to the rank-$r$ root $\alpha _r$. Then we will shown that $\alpha _r$ extends to the  rank-$r$ extremal matroid $\Omega _r$.
\bigskip

\noindent {\bf Claim A.} The rank-$(r-1)$ root  $\alpha _{r-1}$ gives rise to the rank-$r$ root  $\alpha _r$
\bigskip

\noindent {\bf Proof.} By the Strong Splitter Theorem, $M$ must be a cosimple single-element coextension of one of the following matroids: $\alpha _{r-1}$, $\alpha _{r-1, 1}$, $\alpha _{r-1, 2}$, $\alpha _{r-1, 1, 1}$, $\alpha _{r-1, 1, 2}$, or $\alpha _{r-1, 2, 2}$. Moreover, if $M$ is a cosimple single-element coextension of $\alpha _{r-1, 1, 1}$, $\alpha _{r-1, 1, 2}$, or $\alpha _{r-1, 2, 2}$, then $\{b_r, c_r, d_r\}$ is a triad ({\it i.e.}  the matrix is formed with row $[0 0 \dots 0 1 1]$). 
\medskip

First observe that if row $[0 0 \dots 0 1 1]$ is added to $\alpha _{r-1, 1, 1}$ we get precisely $\alpha _r$. Thus if $M$ is a cosimple single-element coextension of $\alpha _{r-1, 1, 1}$ where  $\{b_r, c_r, d_r\}$ is a triad, then  $M=\alpha_r$, which is the  rank-$r$ root matroid. We will show that the other matroids do not have cosimple single-element coextensions in $EX[P_9^*]$.
\medskip

\noindent {\bf Case (i)} By the induction hypothesis $\alpha _{r-1}$ is  formed by adding row $[0 0 \dots 0 1 1]$ to $\alpha _r$ and therefore has no further single-element coextension in $EX[P_9^*]$. Thus, $M$ is not a cosimple single-element coextension of $\alpha _{r-1}$. 
\medskip

\noindent {\bf Case (ii)} Suppose, if possible, $M$ is a cosimple single-element coextension of $\alpha _{r-1, 1}$. Observe that,
$$\alpha_{r-1, 1}\backslash c_{r-1} = \alpha _{r-1}.$$
Therefore, there are no Type(i) rows to be added to $\alpha _{r-1, 1}$. Only Type II and Type III rows can be added to $\alpha _{r-1, 1}$. Type II rows are  the identity rows with a one in the last entry and Type III rows are $b_1, b_2, \dots b_{r-2}, b_{r-1}$ with the last entry switched (put 0 if the last entry is 1 and 1 if it is 0). The superscripts indicate if the last entry is a 1 or a 0. 
Thus, using $\alpha_{r-1, 1}\backslash c_{r-1} = \alpha _{r-2}$ the choices for the last row are:

\medskip
\noindent {\bf Type II rows}:\footnote {Note to referee: All these rows can be written in paragraph format making the paper much shorter. It is displayed like this to help the referee follow the intricate arguments.} 
$$a_1^1=[1 0 0 \dots 0 0 0 0 {\bf 1}]$$ 
$$a_2^1=[0 1 0 \dots 0 0 0 0 {\bf 1}]$$ 
$$.$$
$$.$$
$$c_{r-2}^1=[ 0 0 0 \dots 0 0 1 0 {\bf 1}]$$
$$d_{r-2}^1=[0 0 0 \dots 0 0 0 1 {\bf 1}]$$

\medskip
\noindent {\bf Type III rows}: 
$$b_1^0=[0 1 1 1 1 1 0 \dots 1 0 1 0 {\bf 0}]$$
$$b_2^0=[1 0 1 1 1 1 0 \dots 1 0 1 0 {\bf 0}]$$ 
$$b_3^1=[1 1 0 1 0 0 1 \dots 0 1 0 1 {\bf 1}]$$ 
$$b_4^1=[1 1 1 1 0 0 1 \dots 0 1 0 1 {\bf 1}]$$  
$$b_5^1=[0 0 0 1 1 0 0 \dots 0 1 0 1 {\bf 1}]$$
$$b_6^1=[0 0 0 0 0 1 1 \dots 0 0 0 0 {\bf 1}]$$
$$.$$
$$.$$
$$b_{r-2}^1=[0 0 0 0 0 0 0 \dots  1 1 0 0 {\bf 1}]$$
$$b_{r-1}^1=[0 0 0 0 0 0 0 \dots  0 0 1 1 {\bf 1}]$$

\medskip
\noindent On the other hand, observe that $$\alpha_{r-1, 1}\backslash c_{r-2}/b_{r-1} \cong \alpha _{r-2, 1, 1}.$$
Note the isomorphism instead of inequality. This happens because in $\alpha_{r-1, 1}\backslash c_{r-2}/b_{r-1}$ the last two columns $d_{r-2}$ and $c_{r-1}$ are switched. Otherwise it would be exactly equal to $\alpha _{r-2, 1, 1}$. By the induction hypothesis, $\alpha _{r-2, 1, 1}$ has exactly one cosimple single-element coextension in the class ({\it viz.} the one formed by row $x=[0 0 0 \dots 0 1 1]$). The isomorphism instead of equality is of no consequence since the last two entries in a Type I row are both ones.  Thus we have Type I rows with the a zero or one in the third last entry (the position of the deleted column $c_{r-2}$), Type II rows which are the identity rows with a one in the third last entry, and Type III rows with the third last entry switched.

\medskip
\noindent {\bf Type I rows}: 
$x^0=[0 0 0 \dots 0 {\bf 0} 1 1]$, 
$x^1=[0 0 0 \dots 0 {\bf 1} 1 1]$ 

\medskip
\noindent {\bf Type II rows}: 
$$a_1^1=[1 0 0 \dots 0 0 {\bf 1} 0 0]$$ 
$$a_2^1=[0 1 0 \dots 0 0 {\bf 1} 0 0]$$ 
$$.$$
$$.$$
$$c_{r-3}^1=[ 0 0 0 \dots 0 1 {\bf 1} 0 0]$$
$$d_{r-3}^1=[ 0 0 0 \dots 0 0 {\bf 1} 1 0]$$ 
$$d_{r-2}^1=[0 0 0 \dots 0 0 {\bf 1} 0 1]$$

\medskip
\noindent {\bf Type III rows}: 
$$b_1^0=[0 1 1 1 1 1 1\dots 1 0 {\bf 0} 0 1]$$
$$b_2^0=[1 0 1 1 1 1 1\dots 1 0 {\bf 0} 0 1]$$ 
$$b_3^1=[1 1 0 1 0 0 0 \dots 0 1 {\bf 1} 1 0]$$ 
$$b_4^1=[1 1 1 1 0 0 0 \dots 0 1 {\bf 1} 1 0]$$  
$$b_5^1=[0 0 0 1 1 0 0 \dots 0 0 {\bf 1} 0 0]$$
$$b_6^1=[0 0 0 0 0 1 1 \dots 0 0 {\bf 1} 0 0]$$
$$.$$
$$.$$
$$b_{r-2}^1=[0 0 0 0 0 0 0 \dots 1 1 {\bf 1} 0 0]$$
 
The only common rows are $[0 0 0 \dots 0 0 1 0 1]$, $[0 0 0 \dots 0 0 0 1 1]$ and $[0 0 0 \dots 0 1 1 1]$. Therefore, the only matrices that must be checked explicitly for a $P_9$ minor are the ones formed with the above three rows. These three rank-$r$ matroids have the following three rank-7 minors, respectively, obtained by contracting $\{b_7, \dots , b_{r-1}\}$ and deleting $\{c_5, d_5, c_6, d_6, \dots , c_{r-3}, d_{r-3}\}$.
\bigskip

\tiny
\[ 
M_1=\left[ 
\begin{array}{c|cccccccc}
&    0&1&1&1&1&1&0&1 \\
&    1&0&1&1&1&1&0&1 \\
I_7& 1&1&0&1&0&0&1&0 \\
&    1&1&1&1&0&0&1&0 \\
&    0&0&0&1&1&0&0&0 \\
&    0&0&0&0&0&1&1&0 \\
&    0&0&0&0&0&0&1&1  
\end{array} 
\right ] 
M_2=\left[ 
\begin{array}{c|cccccccc}
&    0&1&1&1&1&1&0&1 \\
&    1&0&1&1&1&1&0&1 \\
I_7& 1&1&0&1&0&0&1&0 \\
&    1&1&1&1&0&0&1&0 \\
&    0&0&0&1&1&0&0&0 \\
&    0&0&0&0&0&1&1&0 \\
&    0&0&0&0&0&1&1&1  
\end{array} 
\right] 
M_3=\left[ 
\begin{array}{c|cccccccc}
&    0&1&1&1&1&1&0&1 \\
&    1&0&1&1&1&1&0&1 \\
I_7& 1&1&0&1&0&0&1&0 \\
&    1&1&1&1&0&0&1&0 \\
&    0&0&0&1&1&0&0&0 \\
&    0&0&0&0&0&1&1&0 \\
&    0&0&0&0&0&1&0&1  
\end{array} 
\right] 
\]
\normalsize
\medskip

\noindent The above matroids have a $P_9^*$-minor. Thus $M$ cannot be a cosimple single-element coextension of $\alpha _{r-1, 1}$.
\medskip

\noindent {\bf Case (iii)} Suppose, if possible, $M$ is a cosimple single-element coextension of $\alpha _{r-1, 2}$. Since
$$\alpha_{r-1, 2}\backslash e_{r-1} = \alpha _{r-1}.$$
there are no Type I rows to be added to $\alpha _{r-1, 1}$ for the same reason as in the previous case. Type II and Type III rows that may be added to $\alpha _{r-1}$ are:
\medskip

\noindent {\bf Type II rows}: 
$$a_1^1=[1 0 0 \dots 0 0 0 0 {\bf 1}]$$  
$$a_2^1=[0 1 0 \dots 0 0 0 0 {\bf 1}]$$ 
$$.$$
$$.$$
$$c_{r-2}^1=[ 0 0 0 \dots 0 0 1 0 {\bf 1}]$$  
$$d_{r-2}^1=[0 0 0 \dots  0 0 0 1 {\bf 1}]$$ 

\medskip
\noindent{\bf Type III rows}: 
$$b_1^0=[0 1 1 1 1 1 0 \dots 1 0 1 0 {\bf 0}]$$
$$b_2^0=[1 0 1 1 1 1 0 \dots 1 0 1 0 {\bf 0}]$$ 
$$b_3^0=[1 1 0 1 0 0 1 \dots 0 1 0 1 {\bf 0}]$$ 
$$b_4^1=[1 1 1 1 0 0 1 \dots 0 1 0 1 {\bf 1}]$$  
$$b_5^1=[0 0 0 1 1 0 0 \dots 0 0 0 0 {\bf 1}]$$
$$b_6^1=[0 0 0 0 0 1 1 \dots 0 0 0 0 {\bf 1}]$$
$$.$$
$$.$$
$$b_{r-2}^1=[0 0 0 0 0 0 0 \dots  1 1 0 0 {\bf 1}]$$
$$b_{r-1}^1=[0 0 0 0 0 0 0 \dots  0 0 1 1 {\bf 1}]$$
\medskip

\noindent Similarly, observe that,
$$\alpha_{r-1, 2}\backslash d_{r-2}/b_{r-1} = \alpha _{r-2, 1, 2}.$$
There are no Type(i) rows to be added to $\alpha _{r-1, 1}$ by the induction hypothesis. Type II and Type III rows that may be added to $\alpha _{r-1}$ are:
\medskip

\noindent{\bf Type II rows}: 
$$a_1^1=[1 0 0 \dots 0 0 0 {\bf 1} 0]$$ 
$$a_2^1=[0 1 0 \dots 0 0 0 {\bf 1} 0]$$ 
$$.$$
$$.$$
$$c_{r-2}^1=[ 0 0 0 \dots 0 0 1 {\bf 1} 0]$$ 
$$e_{r-1}^1=[0 0 0 \dots 0 0 0 0 {\bf 1} 1]$$ 

\medskip
\noindent {\bf Type III rows}: 
$$b_1^0=[0 1 1 1 1 1 0 \dots 1 0 1 {\bf 1} 1]$$
$$b_2^0=[1 0 1 1 1 1 0 \dots 1 0 1 {\bf 1} 1]$$ 
$$b_3^1=[1 1 0 1 0 0 1 \dots 0 1 0 {\bf 0} 0]$$ 
$$b_4^1=[1 1 1 1 0 0 1 \dots 0 1 0 {\bf 0} 0]$$  
$$b_5^1=[0 0 0 1 1 0 0 \dots 0 0 0 {\bf 1} 0]$$
$$b_6^1=[0 0 0 0 0 1 1 \dots 0 0 0 {\bf 1} 0]$$
$$.$$
$$.$$
$$b_{r-2}^1=[0 0 0 0 0 0 0 \dots  1 1 0 0 {\bf 1} 0]$$
\medskip

\noindent The only common row is $[0 0 0 \dots 0 0 0 1 1]$. Therefore, the only matrix that must be checked explicitly for a $P_9$ minor is the following matrix with row  $[0 0 0 \dots 0 0 0 1 1]$. This matrix has the following rank-7 minor obtained by obtained by contracting $\{b_7, \dots , b_{r-1}\}$ and deleting $\{c_5, d_5, c_6, d_6, \dots , c_{r-3}, d_{r-3}\}$.  
\bigskip

\small
\[ 
M_4=\left[ 
\begin{array}{c|cccccccc}
&    0&1&1&1&1&1&0&1 \\
&    1&0&1&1&1&1&0&1 \\
I_7& 1&1&0&1&0&0&1&1 \\
&    1&1&1&1&0&0&1&0 \\
&    0&0&0&1&1&0&0&0 \\
&    0&0&0&0&0&1&1&0 \\
&    0&0&0&0&0&0&1&1  
\end{array} 
\right] 
\]
\normalsize
\medskip

\noindent The above matroid has a $P_9^*$-minor. Thus $M$ cannot be a cosimple single-element coextension of  $\alpha _{r-1, 2}$.
\medskip

\noindent {\bf Case (iv)} Suppose, if possible, $M$ is a cosimple single-element coextension of $\alpha _{r-1, 1, 2}$ or $\alpha _{r-1, 2, 2}$. Then $M$ is formed by adding row $[0 0 0 \dots 0 1 1]$ to $\alpha _{r-1, 1, 2}$ or $\alpha _{r-1, 2, 2}$. The matrices formed in this manner denoted by $\alpha _{r-1, 1, 2}'$ or $\alpha _{r-1, 2, 2}'$ are shown below. They have a $P_9^*$-minor. 

\tiny
\[ 
\alpha _{r-1, 1, 2}'=\left[ 
\begin{array}{c|ccccccccc}
&    0&1&1&1&1&1&0&1&1 \\
&    1&0&1&1&1&1&0&1&1 \\
I_7& 1&1&0&1&0&0&1&0&1 \\
&    1&1&1&1&0&0&1&0&0 \\
&    0&0&0&1&1&0&0&0&0 \\
&    0&0&0&0&0&1&1&0&0 \\
&    0&0&0&0&0&0&0&1&1  
\end{array} 
\right]
\alpha _{r-1, 2, 2}'=\left[ 
\begin{array}{c|ccccccccc}
&    0&1&1&1&1&1&0&1&0 \\
&    1&0&1&1&1&1&0&1&0 \\
I_7& 1&1&0&1&0&0&1&1&1 \\
&    1&1&1&1&0&0&1&0&1 \\
&    0&0&0&1&1&0&0&0&1 \\
&    0&0&0&0&0&1&1&0&0 \\
&    0&0&0&0&0&0&0&1&1  
\end{array} 
\right]  
\]
\normalsize

\bigskip

\noindent {\bf Claim B.} The rank-$r$ root $\alpha _r$ extends to the rank-$r$ extremal matroid $\Omega _r$
\bigskip

\noindent {\bf Proof.} We will prove that the only columns that can be added to $\alpha _r$  are  $c_r, d_r, e_r, f_r, g_{r, 1}, \dots , g_{r, r-4}$. Observe that adding all these columns give $\Omega _r$.

We begin by showing that $\alpha _r$ has two single-element extensions $\alpha _{r, 1}$ and $\alpha _{r, 2}$. Observe that $$\alpha _r / b_r = \alpha _{r-1}.$$ By the induction hypothesis the only columns that can be added to $\alpha _{r-1}$ are $e_{r-1}, f_{r-1}, g_{r, 1}, \dots , g_{r, r-5}$. There are three types of columns that can be added to $\alpha _r$. 

Type I columns are  $\alpha _{r-1}$ are $e_{r-1}, f_{r-1}, g_{r, 1}, \dots , g_{r, r-5}$ with a zero or one in the last entry. Type II and III columns are the columns of $\alpha _{r-1}$ with a zero or one in the last entry. The columns are listed below with the last entry highlighted:

\medskip
\noindent {\bf Type I columns:} 
$$e_{r-1}^0=[1 1 1 0 0 0 0  \dots 0 0  {\bf 0}]^T$$
$$e_{r-1}^1=[1 1 1 0 0 0 0  \dots 1 0  {\bf 1}]^T$$ 
$$f_{r-1}^0=[0 0 1 1 1 0 0  \dots 0 1  {\bf 0}]^T$$ 
$$e_{r-1}^1=[0 0 1 1 1 0 0  \dots 0 1  {\bf 1}]^T$$  
$$g_{r-1, 1}^0=[1 1 1 1 0 0 0 \dots 0 1  {\bf 0}]^T$$
$$g_{r-1, 1}^1=[1 1 1 1 0 0 0 \dots 0 0  {\bf 1}]^T$$
$$g_{r-1, 2}^0=[1 1 1 1 0 1 0 \dots 0 0  {\bf 0}]^T$$
$$g_{r-1, 2}^1=[1 1 1 1 0 1 0 \dots 0 0  {\bf 1}]^T$$
$$.$$
$$.$$
$$g_{r-1, r-5}^0=[1 1 1 1 0 0 0 \dots 0 1  {\bf 0}]^T$$
$$g_{r-1, r-5}^1=[1 1 1 1 0 0 0 \dots 0 1  {\bf 1}]^T$$
\medskip

\noindent {\bf Type II Columns}: 
$$b_1^1=[1 0 0 \dots 0 0 {\bf 1}]^T$$
$$b_2^1=[0 1 0 \dots 0 0 {\bf 1}]^T$$ 
$$.$$
$$.$$
$$b_{r-1}^1=[ 0 0 0 \dots 0 1 {\bf 1}]^T$$

\medskip
\noindent{\bf Type III Columns}: 
$$a_1^0=[0 1 1 1 0 0 0 \dots 0 0  {\bf 1}]^T$$
$$a_2^0=[1 0 1 1 0 0 0 \dots 0 0  {\bf 1}]^T$$ 
$$a_3^1=[1 1 0 1 0 0 0 \dots 0 0  {\bf 1}]^T$$ 
$$a_4^1=[1 1 1 1 1 0 0 \dots 0 0  {\bf 1}]^T$$  
$$a_5^1=[1 1 0 0 1 0 0 \dots 0 0  {\bf 1}]^T$$
$$c_5^1=[1 1 0 0 0 1 0 \dots 0 0  {\bf 1}]^T$$
$$d_5^1=[0 0 1 1 0 1 0 \dots 0 0  {\bf 1}]^T$$
$$c_6^1=[1 1 0 0 0 0 1 \dots 0 0  {\bf 1}]^T$$
$$d_6^1=[0 0 1 1 0 0 1 \dots 0 0  {\bf 1}]^T$$
$$.$$
$$.$$
$$c_{r-2}^1=[1 1 0 0 0 0 0 \dots  0 1 {\bf 1}]^T$$
$$d_{r-2}^1=[0 0 1 1 0 0 0 \dots  0 1 {\bf 1}]^T$$
$$c_{r-1}^0=[1 1 0 0 0 0 0 \dots  0 0 {\bf 0}]^T$$
$$d_{r-1}^0=[0 0 1 1 0 0 0 \dots  0 0 {\bf 0}]^T$$
\medskip

\noindent Observe that, $c_{r-1}^0=c_r$, $d_{r-1}^0=d_r$  $e_{r-1}^0=e_r$, $f_{r-1}^0=f_r$, 
$g_{r-1, 1}^0=g_{r, 1}, g_{r-1, 2}^0=g_{r, 2}, g_{r-1, r-5}^0=g_{r, r-5}$, and $g_{r-1, 1}^1=g_{r, r-4}$. We must show that the matrices obtained by adding the other columns have $P_9^*$-minor. Following the notation in [\ref{Oxley2012}] we will denote a single-element extension of $M$ as $M+e$.

Consider the Type I columns $e_{r-1}^1$ and $f_{r-1}^1$. Writing out the matrices it is easy to see that
$$(\alpha _r +e_{r-1}^1) / \{b_6, b_7, \dots b_{r-1}\} \backslash \{c_5, d_5, \dots c_{r-2}, d_{r-2}\}=\alpha _6 + e_6^1$$
and
$$(\alpha _r +f_{r-1}^1) / \{b_6, b_7, \dots b_{r-1}\} \backslash \{c_5, d_5, \dots c_{r-2}, d_{r-2}\}=\alpha _6 + f_6^1.$$ 
The matrices $\alpha _6 + e_6^1$ and $\alpha _6 + f_6^1$ are shown below. Both have a $P_9^*$-minor.

\small
\[ 
\alpha _6 + e_6^1=\left[ 
\begin{array}{c|cccccccc}
&    0&1&1&1&1&1&0&1 \\
&    1&0&1&1&1&1&0&1 \\
I_6& 1&1&0&1&0&0&1&1 \\
&    1&1&1&1&0&0&1&0 \\
&    0&0&0&1&1&0&0&0 \\
&    0&0&0&0&0&1&1&1 \\  
\end{array} 
\right]
\alpha _6 + f_6^1=\left[ 
\begin{array}{c|cccccccc}
&    0&1&1&1&1&1&0&0 \\
&    1&0&1&1&1&1&0&0 \\
I_6& 1&1&0&1&0&0&1&1 \\
&    1&1&1&1&0&0&1&1 \\
&    0&0&0&1&1&0&0&1 \\
&    0&0&0&0&0&1&1&1 \\  
\end{array} 
\right]  
\]
\normalsize

\noindent Consider the Type I columns $g_{r-1, 2}^1, g_{r-1, 3}^1, \dots g_{r-1, r-5}^1$. For $2 \le k \le r-5$, the matrix $\alpha _r +g_{r-1, k}^1$ has as minor $\alpha _7 + g_{6, 1}^1$ obtained by contracting all  columns $\{b_6, \dots , b_{r-1}\}$ except $b_{k+4}$ and deleting all columns $\{c_5, d_5, \dots , c_{r-2}, d_{r-2}\}$ except $c_{k+3}$ and $d_{k+3}$.
The matrix $\alpha _7 + g_{6, 1}^1$ is shown below. It has a $P_9^*$-minor.
\small
\[ 
\alpha _7 + g_{6, 1}^1=\left[ 
\begin{array}{c|cccccccccc}
&    0&1&1&1&1&1&0&1&0&1 \\
&    1&0&1&1&1&1&0&1&0&1 \\
I_7& 1&1&0&1&0&0&1&0&1&1 \\
&    1&1&1&1&0&0&1&0&1&1 \\
&    0&0&0&1&1&0&0&0&0&0 \\
&    0&0&0&0&0&1&1&0&0&1 \\
&    0&0&0&0&0&0&0&1&1&1  
\end{array} 
\right]  
\]
\normalsize

Consider the Type II columns $b_1^1, \dots b_{r-1}^1$. Writing out the matrices we see that for $1\le k \le 5$
$$(\alpha _r +b_k^1) / \{b_6, b_7, \dots b_{r-1}\} \backslash \{c_5, d_5, \dots c_{r-2}, d_{r-2}\}=\alpha _6 + b_k^1.$$ 
These matrices are shown below. They have a $P_9^*$-minor.

\tiny
\[ 
\alpha _6 + b_1^1=\left[ 
\begin{array}{c|cccccccc}
&    0&1&1&1&1&1&0&1 \\
&    1&0&1&1&1&1&0&0 \\
I_6& 1&1&0&1&0&0&1&0 \\
&    1&1&1&1&0&0&1&0 \\
&    0&0&0&1&1&0&0&0 \\
&    0&0&0&0&0&1&1&1 \\  
\end{array} 
\right]
\alpha _6 + b_2^1=\left[ 
\begin{array}{c|cccccccc}
&    0&1&1&1&1&1&0&0 \\
&    1&0&1&1&1&1&0&1 \\
I_6& 1&1&0&1&0&0&1&0 \\
&    1&1&1&1&0&0&1&0 \\
&    0&0&0&1&1&0&0&0 \\
&    0&0&0&0&0&1&1&1 \\  
\end{array} 
\right]
\]
\normalsize

\tiny
\[  
\alpha _6 + b_3^1=\left[ 
\begin{array}{c|cccccccc}
&    0&1&1&1&1&1&0&0 \\
&    1&0&1&1&1&1&0&0 \\
I_6& 1&1&0&1&0&0&1&1 \\
&    1&1&1&1&0&0&1&0 \\
&    0&0&0&1&1&0&0&0 \\
&    0&0&0&0&0&1&1&1 \\  
\end{array} 
\right]  
\alpha _6 + b_4^1=\left[ 
\begin{array}{c|cccccccc}
&    0&1&1&1&1&1&0&0 \\
&    1&0&1&1&1&1&0&0 \\
I_6& 1&1&0&1&0&0&1&0 \\
&    1&1&1&1&0&0&1&1 \\
&    0&0&0&1&1&0&0&0 \\
&    0&0&0&0&0&1&1&1 \\  
\end{array} 
\right]
\]
\normalsize

\tiny
\[
\alpha _6 + b_5^1=\left[ 
\begin{array}{c|cccccccc}
&    0&1&1&1&1&1&0&0 \\
&    1&0&1&1&1&1&0&0 \\
I_6& 1&1&0&1&0&0&1&0 \\
&    1&1&1&1&0&0&1&0 \\
&    0&0&0&1&1&0&0&1 \\
&    0&0&0&0&0&1&1&1 \\  
\end{array} 
\right] 
\]
\normalsize

For $6\le k \le r-1$, the matroid $\alpha _r + b_k$ has minor $\alpha _7 + b_6^1$ obtained by contracting all columns $\{b_6, \dots , b_{r-1}\}$ except $b_k$ and deleting all  columns $\{c_5, d_5, \dots , c_{r-2},d_{r-2}$\} except $c_{r-2}$ and  $d_{r-2}$.
The matrix $\alpha _7 + b_6^1$ is shown below. It has a $P_9^*$-minor.

\small
\[ 
\alpha _7 + b_6^1=\left[ 
\begin{array}{c|cccccccccc}
&    0&1&1&1&1&1&0&1&0&0 \\
&    1&0&1&1&1&1&0&1&0&0 \\
I_7& 1&1&0&1&0&0&1&0&1&0 \\
&    1&1&1&1&0&0&1&0&1&0 \\
&    0&0&0&1&1&0&0&0&0&0 \\
&    0&0&0&0&0&1&1&0&0&1 \\
&    0&0&0&0&0&0&0&1&1&1  
\end{array} 
\right]  
\]
\normalsize

Consider the Type III columns $a_1^1, a_2^1, a_3^1, a_4^1, a_5^1, c_5^1, d_5^1, c_6^1, d_6^1, \dots c_{r-1}^0, d_{r-1}^0$. Writing out the matrices we see that for $1\le k \le 5$
$$\alpha _r +a_k^1 / \{b_6, b_7, \dots b_{r-1}\} \backslash \{c_5, d_5, \dots c_{r-2}, d_{r-2}\}=\alpha _6 + a_k^1$$
These matrices are shown below. They have a $P_9^*$-minor.

\tiny
\[ 
\alpha _6 + a_1^1=\left[ 
\begin{array}{c|cccccccc}
&    0&1&1&1&1&1&0&0 \\
&    1&0&1&1&1&1&0&1 \\
I_6& 1&1&0&1&0&0&1&1 \\
&    1&1&1&1&0&0&1&1 \\
&    0&0&0&1&1&0&0&0 \\
&    0&0&0&0&0&1&1&1 \\  
\end{array} 
\right]
\alpha _6 + a_2^1=\left[ 
\begin{array}{c|cccccccc}
&    0&1&1&1&1&1&0&1 \\
&    1&0&1&1&1&1&0&0 \\
I_6& 1&1&0&1&0&0&1&1 \\
&    1&1&1&1&0&0&1&1 \\
&    0&0&0&1&1&0&0&0 \\
&    0&0&0&0&0&1&1&1 \\  
\end{array} 
\right]  
\]
\normalsize

\tiny
\[
\alpha _6 + a_3^1=\left[ 
\begin{array}{c|cccccccc}
&    0&1&1&1&1&1&0&1 \\
&    1&0&1&1&1&1&0&1 \\
I_6& 1&1&0&1&0&0&1&0 \\
&    1&1&1&1&0&0&1&1 \\
&    0&0&0&1&1&0&0&0 \\
&    0&0&0&0&0&1&1&1 \\  
\end{array} 
\right]  
\alpha _6 + a_4^1=\left[ 
\begin{array}{c|cccccccc}
&    0&1&1&1&1&1&0&1 \\
&    1&0&1&1&1&1&0&1 \\
I_6& 1&1&0&1&0&0&1&1 \\
&    1&1&1&1&0&0&1&1 \\
&    0&0&0&1&1&0&0&1 \\
&    0&0&0&0&0&1&1&1 \\  
\end{array} 
\right]
\]
\normalsize

For $6\le k \le r-2$. the matrix $\alpha _r +c_k^1$ has as minor $\alpha _7 + c_5^1$ obtained by contacting columns  $\{b_6, \dots , b_{r-1}\}$ except $b_{k+1}$ and delting columns $\{c_5, d_5, \dots c_{r-2}, d_{r-2}\}$ except $c_k$ and $d_k$. Similary, $\alpha _r +d_k^1$ has as minor $\alpha _7 + d_5^1$. The matrices $\alpha _7 + c_5^1$ and  $\alpha _7 + d_5^1$ are shown below. They have a $P_9^*$-minor.
\small
\[ 
\alpha _7 + c_5^1=\left[ 
\begin{array}{c|cccccccccc}
&    0&1&1&1&1&1&0&1&0&1 \\
&    1&0&1&1&1&1&0&1&0&1 \\
I_7& 1&1&0&1&0&0&1&0&1&0 \\
&    1&1&1&1&0&0&1&0&1&0 \\
&    0&0&0&1&1&0&0&0&0&0 \\
&    0&0&0&0&0&1&1&0&0&1 \\
&    0&0&0&0&0&0&0&1&1&1  
\end{array} 
\right]  
\alpha _7 + d_5^1=\left[ 
\begin{array}{c|cccccccccc}
&    0&1&1&1&1&1&0&1&0&0 \\
&    1&0&1&1&1&1&0&1&0&0 \\
I_7& 1&1&0&1&0&0&1&0&1&1 \\
&    1&1&1&1&0&0&1&0&1&1 \\
&    0&0&0&1&1&0&0&0&0&0 \\
&    0&0&0&0&0&1&1&0&0&1 \\
&    0&0&0&0&0&0&0&1&1&1  
\end{array} 
\right]  
\]
\normalsize
\noindent This completes the proof of Lemma 2.4 that $\alpha _r$ extends to $\Omega _r$. $\qed$
\bigskip
 
Returning to the proof of Theorem 3.1, it is easy to see that $\alpha_r$ has two non-isomorphic single-element extensions $\alpha_{r, 1}$ formed by adding columns $c_r$, $d_r$ or $g_{r,1}$, and the remaining columns give $\alpha_{r, 2}$. The matroid $\alpha_{r, 1}$ has two non-isomorphic single-element extension, the notable one that gives rise to the rank-(r+1) root matroid is $\alpha_{r, 1, 1}$ formed by adding $c_r$ and $d_r$ to $\alpha _r$. $\qed$

\bigskip

We end this paper by showing how Theorem 1.1 follows from Theorem 3.1.
\bigskip

\noindent {\bf Proof of Theorem 1.1.} Each of the matroids and families of matroids listed in Theorem 3.1 is an extremal matroid for $EX[P_9^*]$. Theorem 3.1 implies that the growth rate of $Z_r$ and $\Omega _r$ are the only relevant ones to determine the growth rate of the class $EX[P_9^*]$. The size of $Z_r$ as mentioned earlier is $2r+1$ [\ref{Oxley1987}]. Proposition 2.2 implies that the size of $\Omega _r$ is $4r-5$. Thus the growth-rate function of the non-regular matroids with no $P_9^*$-minor is $4r-5$. It follows that the non-regular members of $EX[P_9^*]$ have linear growth rate. The growth rate function for regular matroids is $\frac {r(r+1)}{2}$ with this bound being attained if and only if $M\cong K_{r+1}$. Hence Theorem 1.1 is proved. $\qed$
\bigskip

%%%%%%%%%%%%%%%%%%%%%%%%%%%%%%%%%%%%%%%%%%%%
\bigskip

\noindent {\bf References}

\bigskip
\begin{enumerate}

\item  \label{Dirac1963} G. A. Dirac (1963). Some results concerning the structure of graphs, {\it Canad. Math. Bull.} {\bf 6}, 183-210.

\item  \label{KinganLemos2002}S. R. Kingan  and M. Lemos (2002). Almost-graphic matroids, {Advances in Applied Mathematics}, {\bf 28}, 438 - 477.

\item  \label{KinganLemos2014}	S. R. Kingan  and M. Lemos  (2014). Strong Splitter Theorem, {\it Annals of Combinatorics} {\bf 18-1}, 111-116.

\item  \label{KinganSubmitted}	S. R. Kingan (submitted). Finding extremal matroids in excluded minor classes.

\item \label{KungMayhewPivottoRoylesubmitted} J. Kung, D. Mayhew, I. Pivotto, and G. Royle (2014). Maximum size binary matroids with no AG(3,2)-minor are graphic, {\it SIAM J. Discrete Math.} {\bf 28(3)}, 1559–1577. 

\item \label{Mader1967} W. Mader (1967). Homomorphieeigenschaften und mittlere kantendichte von graphen. {\it Mathematische Annalen}, {\bf  174}, 265–268.  

\item \label{Oxley1987} J. G. Oxley (1987). The binary matroids with no 4-wheel minor, {\it Trans. Amer. Math. Soc.}  {\bf 154}, 63-75.

\item \label{Oxley2012} J. G. Oxley (2012). {\it Matroid Theory}, Second Edition, Oxford University Press, New York. 

\end{enumerate}
%%%%%%%%%%%%%%%%%%%%%%%%%%%%%%%%%%%%%%%%%%%%%%%%%%%

\vfill
\eject
{\bf Appendix}

\tiny
 \begin{center}
\begin{tabular}{|c|p{15em}|c|}
\hline
\bf{Matroid}& \bf{Extension Columns} & {\bf Name}   \\ \hline  
$P_9$ &  $\bf [1 1 1 0]$  & $D_1$   \\ \hline
&  $\bf [1 0 0 1]$ $[0 1 0 1]$ $[0 1 1 0]$, $[1 0 1 0]$  & $D_2$               \\ \hline
& $\bf [0 0 1 1]$    &   $D_3$          \\  \hline \hline
$D_1$ & $[0 1 0 1]$ $[0 1 1 0]$ $[1 0 0 1]$ $[1 0 1 0]$ & $X_1$ \\  \hline  
& $\bf [0 0 1 1]$ & $X_2$   \\  \hline 
$D_2$ & $\bf [1 0 1 0]$ $[1 1 1 0]$  & $X_1$   \\  \hline
& $\bf [0 0 1 1 ]$ $[0 1 0 1]$ $[0 1 1 0]$ & $X_3$   \\  \hline 
$D_3$ & $[1 1 1 0]$  & $X_2$   \\  \hline 
& $ [0 1 0 1]$ $[0 1 1 0]$ $[1 0 0 1]$ $[1 0 1 0]$ & $X_3$ \\ \hline \hline 
\end{tabular}
 \end{center}
 
 \begin{center} Table 1a: Simple single-element extensions of $P_9$ \end{center} 

 \begin{center}
\begin{tabular}{|p{28em}|c|c}
\hline
\bf{Coextension Rows} & {\bf Name}   & {\bf $P_9^*$-minor}  \\   \hline \hline
$[1 1 0 0 0]$ $[1 1 1 1 1]$    & $E_1$   & Yes \\ \hline
$[1 1 0 1 1]$ $[1 1 1 0 0]$    & $E_2$ & Yes \\  \hline
$[1 1 0 0 1]$ $[1 1 1 0 1]$    & $E_3$ & Yes \\  \hline
$[0 1 0 0 1]$ $[0 1 0 1 0]$ $[0 1 1 0 1]$ $[0 1 1 1 0]$ $[1 0 0 0 1]$ $[1 0 0 1 0]$ $[1 0 1 0 1]$ $[1 0 1 1 0]$    &  $E_4$  & Yes\\  \hline
$[0 1 0 1 1]$ $[0 1 1 0 0]$ $[1 0 0 1 1]$ $[1 0 1 0 0]$   &  $E_5$ & Yes \\  \hline
$[0 0 1 0 1]$ $[0 0 1 1 0]$    &  $E_6$  & Yes\\  \hline
$[0 0 1 1 1]$                  &  $E_6^*$ & Yes \\  \hline
$[0 0 0 1 1]$      & $E_7$ ($\bf \alpha_5$) & No \\  \hline
\end{tabular}
 \end{center}
 
 \begin{center} Table 1b: Cosimple single-element coextensions of $P_9$ \end{center}

 \begin{center}
\begin{tabular}{|p{20em}|c|c|}
\hline
 \bf{Extension Columns} & {\bf Name} & {\bf  $P_9^*$-minor}  \\  \hline \hline
  $[0 0 0 1 1]$ $[0 0 1 0 1]$ $[1 1 1 0 1]$   & Ext 1 &  Yes   \\  \hline

 $[0 0 1 1 0]$ $[1 1 0 0 0]$ $[1 1 1 1 0]$  & Ext 2  ($\bf \alpha_{5, 1}$)  & No   \\  \hline

 $[0 0 1 1 1]$ $[1 1 1 0 0]$   &   Ext 3 ($\bf \alpha_{5, 2}$)  & No   \\  \hline

 $[0 1 0 0 1]$ $[0 1 0 1 0]$ $[0 1 1 0 0]$ $[0 1 1 1 1]$  $[1 0 0 0 1]$ 
$[1 0 0 1 0]$ $[1 0 1 0 0]$ $[1 0 1 1 1]$ & Ext 4   & Yes   \\  \hline

 $[0 1 0 1 1]$ $[0 1 1 0 1]$ $[1 0 0 1 1]$ $[1 0 1 0 1]$ & Ext 5 ($\bf \alpha_{5, 3}$) & No    \\  \hline

 $[1 1 0 1 1]$   & Ext 6 & Yes   \\  \hline  \hline
\end{tabular}
 \end{center}
\begin{center}  Table 2a: Simple single-element extensions of $\alpha_5$   \end{center} 
\bigskip

\begin{center}
\begin{tabular}{|p{20em}|c|c|}
\hline
\bf{Coextension Rows} & {\bf Name}  & {\bf  $P_9^*$-minor}    \\ \hline \hline

$[0 0 1 0 1]$ $[0 0 1 1 0]$ $[1 1 0 0 1]$ $[1 1 0 1 1]$ $[1 1 1 0 0]$ $[1 1 1 0 1]$  & Coext 1 & Yes  \\  \hline

$[0 0 1 1 1]$ $[1 1 0 0 0]$ $[1 1 1 1 1]$ & Coext 2  & Yes    \\  \hline

$[0 1 0 0 1]$ $[0 1 0 1 0]$ $[0 1 0 1 1]$ $[0 1 1 0 0]$ $[0 1 1 0 1]$ $[0 1 1 1 0]$ $[1 0 0 0 1]$ $[1 0 0 1 0]$ $[1 0 0 1 1]$ $[1 0 1 0 0]$ $[1 0 1 0 1]$ $[1 0 1 1 0]$   &   Coext 3   & Yes  \\  \hline \hline
  
\end{tabular}
\end{center}
 
 \begin{center}   Table 2b: Cosimple single-element coextensions of $\alpha_5$ \end{center} 

 \begin{center}
\begin{tabular}{|c|p{30em}|c|c|}
\hline
\bf {Name} & \bf{Coextension Rows} & {\bf Name}  & {\bf  $P_9^*$-minor}    \\ \hline \hline

$\alpha_{5,1}$ & $\bf [ 0 0 0 0 1 1 ]$
$\bf [ 0 0 0 1 0 1 ]$
$\bf [ 1 1 0 1 0 0 ]$
$\bf[ 1 1 1 1 0 0 ]$
& Coext 1 & Yes  \\  \hline

&$\bf [ 0 0 0 1 1 1 ]$
$\bf [ 0 0 1 0 0 1]$
& Coext 2  & Yes    \\  \hline

&$\bf [ 0 1 0 0 0 1 ]$
$\bf [ 0 1 1 1 1 1]$
$\bf [ 1 0 0 0 0 1 ]$
$\bf [ 1 0 1 1 1 1 ]$
& Coext 7    & Yes  \\  \hline
\hline 

$\alpha_{5,2}$ & $\bf [ 0 0 0 0 1 1 ]$
$\bf [ 0 0 0 1 0 1 ]$
$\bf [ 0 0 0 1 1 0 ]$
$[ 0 0 1 1 1 0 ]$
$[ 1 1 0 0 0 0 ]$
$[ 1 1 1 1 1 0 ]$
& Coext 1 & Yes  \\  \hline

&$\bf [ 0 0 1 0 0 1 ]$
$[ 0 0 1 0 1 0 ]$
$[ 0 0 1 1 0 0 ]$
$[ 0 0 1 1 1 1 ]$
$[ 0 1 0 0 1 1 ]$
$[ 0 1 0 1 0 1 ]$
$[ 0 1 0 1 1 0 ]$
$[ 0 1 1 0 0 0 ]$
$[ 0 1 1 0 1 1 ]$
$[ 0 1 1 1 0 1 ]$
$[ 1 0 0 0 1 1 ]$
$[ 1 0 0 1 0 1 ]$
$[ 1 0 0 1 1 0 ]$
$[ 1 0 1 0 0 0 ]$
$[ 1 0 1 0 1 1 ]$
$[ 1 0 1 1 0 1 ]$
$[ 1 1 0 0 0 1 ]$
$[ 1 1 0 0 1 0 ]$
$\bf [ 1 1 0 1 0 0 ]$
$[ 1 1 0 1 1 1 ]$
$[ 1 1 1 0 0 1 ]$
$[ 1 1 1 0 1 0 ]$
$\bf [ 1 1 1 1 0 0 ]$
$[ 1 1 1 1 1 1 ]$
& Coext 2  & Yes    \\  \hline

&$\bf [ 0 1 0 0 0 1 ]$
$[ 0 1 0 0 1 0 ]$
$[ 0 1 0 1 0 0 ]$
$[ 0 1 0 1 1 1 ]$
$[ 0 1 1 0 0 1 ]$
$[ 0 1 1 0 1 0 ]$
$[ 0 1 1 1 0 0 ]$
$\bf [ 0 1 1 1 1 1 ]$
$\bf [ 1 0 0 0 0 1 ]$
$[ 1 0 0 0 1 0 ]$
$[ 1 0 0 1 0 0 ]$
$[ 1 0 0 1 1 1 ]$
$[ 1 0 1 0 0 1 ]$
$[ 1 0 1 0 1 0 ]$
$[ 1 0 1 1 0 0 ]$
$\bf [ 1 0 1 1 1 1 ]$
& Coext 4    & Yes  \\  \hline
\hline  

$\alpha_{5,3}$ & $\bf [ 0 0 0 0 1 1 ]$
$\bf [ 0 0 0 1 0 1 ]$
$\bf [ 0 0 0 1 1 0 ]$
$\bf [ 0 0 1 0 0 1 ]$
$[ 0 0 1 1 1 1 ]$
$[ 0 1 0 0 1 0 ]$
$[ 0 1 0 1 0 0 ]$
$[ 1 0 0 1 0 0 ]$
$[ 1 0 0 1 1 1 ]$
$[ 1 0 1 0 0 0 ]$
$[ 1 0 1 1 0 1 ]$
$[ 1 1 0 0 0 0 ]$
$\bf [ 1 1 0 1 0 1 ]$
$\bf [ 1 1 1 1 0 0 ]$
$[ 1 1 1 1 1 1 ]$
& Coext 1 & Yes  \\  \hline

&$[ 0 0 1 0 1 0 ]$
$[ 0 0 1 1 0 0 ]$
$\bf [ 0 1 0 0 0 1 ]$
$[ 0 1 0 1 1 1 ]$
$[ 0 1 1 0 0 0 ]$
$[ 0 1 1 0 1 1 ]$
$[ 0 1 1 1 0 1 ]$
$\bf [ 1 0 0 0 0 1 ]$
$[ 1 0 0 0 1 0 ]$
$[ 1 0 1 0 1 1 ]$
$\bf [ 1 0 1 1 1 0 ]$
$[ 1 1 0 0 1 1 ]$
$[ 1 1 0 1 1 0 ]$
$[ 1 1 1 0 0 1 ]$
$[ 1 1 1 0 1 0 ]$
& Coext 2  & Yes    \\  \hline

&$\bf [ 0 0 1 1 1 0 ]$
$[ 0 1 0 1 1 0 ]$
$[ 0 1 1 0 0 1 ]$
$[ 0 1 1 0 1 0 ]$
$[ 0 1 1 1 0 0 ]$
$\bf [ 0 1 1 1 1 1 ]$
$[ 1 0 0 0 1 1 ]$
$[ 1 0 1 0 1 0 ]$
$[ 1 1 0 0 0 1 ]$
$[ 1 1 1 1 1 0 ]$
& Coext 4    & Yes  \\  
\hline \hline  
\end{tabular}
\end{center}
 
 \begin{center}   Table 3: Selected cosimple single-element coextensions of $\alpha_{5,1}$, $\alpha_{5,2}$, and $\alpha_{5,3}$ \end{center} 
\normalsize

\tiny
 \begin{center}
\begin{tabular}{|c|p{25em}|c|c|}
\hline
Matroid &\bf{Coextension Rows} & {\bf Name} & Relevant minors  \\      \hline \hline

$D_1$ & $\bf [0 0 0 0 1 1]$ $\bf [0 0 0 1 0 1]$ $[0 0 1 0 1 0]$ $[0 0 1 1 0 0]$ $[0 1 0 0 1 0]$ $[0 1 0 1 0 0]$ $[0 1 1 0 1 1]$ $[0 1 1 1 0 1]$ $[1 0 0 0 1 0]$ $[1 0 0 1 0 0]$ $[1 0 1 0 1 1]$ $[1 0 1 1 0 1]$ $[1 1 0 0 0 1]$ $[1 1 0 1 1 1]$ $[1 1 1 0 0 0]$ $[1 1 1 1 1 0]$   &    $A_{11}$             & $E_1$, $E_2$, $E_3$, $E_4$, $E_6^*$     \\  \hline

& $\bf [0 0 0 1 1 0]$     &  $A_{30}$             &   $E_7$   \\  \hline

& $\bf [0 0 0 1 1 1]$ $[0 0 1 1 1 0]$ $[0 1 0 1 1 0]$ $[0 1 1 0 0 1]$ $[1 0 0 1 1 0]$ $[1 0 1 0 0 1]$ $[1 1 0 0 1 1]$ $[1 1 1 0 1 0]$   & $A_{18}$     &  $E_3$, $E_5$, $E_6^*$, $E_7$  \\  \hline

& $\bf [0 0 1 0 0 1]$ $[0 0 1 1 1 1]$  &    $A_{14}$           &  $E_2$, $E_6^*$    \\  \hline

& $[0 0 1 0 1 1]$ $[0 0 1 1 0 1]$  &     $A_{13}$        &  $E_2$, $E_3$ $E_5$    \\  \hline

& $\bf [0 1 0 0 0 1]$ $[0 1 0 0 1 1]$  $[0 1 0 1 0 1]$ $[0 1 0 1 1 1]$ $[0 1 1 0 0 0]$ $[0 1 1 0 1 0]$ $[0 1 1 1 0 0]$ $\bf [0 1 1 1 1 0]$ $\bf [1 0 0 0  1]$ $[1 0 0 0 1 1]$ $[1 0 0 1 0 1]$ $[1 0 0 1 1 1]$ $[1 0 1 0 0 0]$ $[1 0 1 0 1 0]$ $[1 0 1 1 0 0 ]$ $\bf [1 0 1 1 1 0]$ &  $A_{15}$ &  $E_2$, $E_6^*$   \\  \hline

& $[1 1 0 0 0 0]$ $\bf [1 1 0 1 0 0]$  $\bf [1 1 1 1 0 1]$ $[1 1 1 1 1 1]$   & $A_{10}$   &      $E_1$, $E_2$ \\  \hline

& $[1 1 0 0 1 0]$ $[1 1 0 1 1 0]$  $[1 1 1 0 0 1]$ $[1 1 1 0 1 1]$    &  $A_{17}$    &  $E_2$, $E_3$  \\  \hline \hline

$D_2$ & $\bf [0 0 0 0 1 1]$ $\bf [0 0 0 1 0 1]$ $\bf [0 0 0 1 1 0]$ $[0 0 1 1 1 1]$  $[1 0 0 1 1 1]$ $[1 0 1 0 0 0]$   &       $A_{26}$ ${\bf A}$          & $E_5$, $E_6^*$, $E_7$   \\  \hline

& $\bf [0 0 0 1 1 1]$     &     $A_{31}$ ${\bf Z}$            &  $E_7$, $R_{10}$   \\  \hline

& $\bf [0 0 1 0 0 1]$ $[0 1 0 1 0 0]$ $[0 1 1 1 0 1]$    & $A_{23}$  & $E_4$, $E_5$  \\  \hline

& $[0 0 1 0 1 0]$ $[0 0 1 1 0 0]$ $\bf [0 1 0 0 0 1]$ $[0 1 0 0 1 0]$ $[0 1 1 0 1 1]$ $\bf [0 1 1 1 1 0]$ &     $A_{20}$            & $E_4$, $E_6$    \\  \hline

& $[0 0 1 0 1 1]$ $[0 0 1 1 0 1]$  $[0 1 0 1 0 1]$ $[0 1 0 1 1 0]$ $[0 1 1 0 0 1]$ $[0 1 1 1 0 0]$ &  $A_{21}$   & $E_4$, $E_5$    \\  \hline

& $[0 0 1 1 1 0]$ $[0 1 0 0 1 1]$  $[0 1 1 0 1 0]$  &    $A_{24}$  & $E_4$    \\  \hline

& $\bf [1 0 0 0 0 1]$ $[1 0 1 0 0 0]$  $[1 0 1 0 1 1]$ $[1 0 1 1 0 1]$ $[1 1 0 1 1 0]$ $[1 1 1 0 0 1]$  &  $A_{15}$   & $E_2$, $E_5$    \\  \hline

& $[1 0 0 0 1 0]$ $[1 0 0 1 0 0]$  $[1 1 0 0 0 0]$ $\bf [1 1 0 1 0 1]$ $\bf [1 1 1 1 0 0]$ $[1 1 1 1 1 1]$  &  $A_{6}$   & $E_1$, $E_4$    \\  \hline

& $[1 0 0 0 1 1]$ $[1 0 0 1 0 1]$  $[1 1 0 0 1 0]$ $[1 1 0 1 1 1]$ $[1 1 1 0 0 0]$ $[1 1 1 0 1 1]$ &    $A_{16}$  & $E_2$, $E_3$, $E_4$, $E_6^*$    \\  \hline

& $[1 0 0 1 1 0]$ $[1 0 1 0 1 0]$  $[1 0 1 1 0 0]$ $\bf [1 0 1 1 1 1]$ $[1 1 0 0 0 1]$ $[1 1 1 1 1 0]$ &    $A_7$  & $E_4$, $E_5$   \\  \hline

& $[1 0 0 1 1 1]$ $[1 1 0 0 1 1]$  $[1 1 1 0 1 0]$   &    $A_{18}$ ${\bf C}$ & $E_3$, $E_5$, $E_6^*$, $E_7$    \\  \hline

& $[1 0 1 0 0 1]$  &    $A_{27}$ $\bf B$  & $E_5$    \\  \hline \hline

$D_3$ & $\bf [0 0 0 0 1 1]$ $\bf [0 0 0 1 0 1]$ $\bf [0 0 0 1 1 0]$     &    $A_{29}$ &  $E_7$  \\  \hline

& $\bf [0 0 0 1 1 1]$     &  $A_{30}$  &   $E_7$   \\  \hline

& $\bf [0 0 1 0 0 1]$ $[1 1 0 1 1 1]$ $[1 1 1 0 0 1]$    &      $A_{14}$  & $E_2$, $E_6^*$     \\  \hline

& $[0 0 1 0 1 0]$ $[0 0 1 1 0 0]$ $[1 1 0 0 0 1]$ $\bf [1 1 0 1 0 0]$ $\bf [1 1 1 1 0 0]$ $[1 1 1 1 1 1]$  &                $A_{4}$ & $E_1$, $E_6$    \\  \hline

& $[0 0 1 0 1 1]$ $[0 0 1 1 0 1]$  $[1 1 0 0 1 1]$ $[1 1 0 1 1 0]$ $[1 1 1 0 0 0]$ $[1 1 1 0 1 1]$ &   $A_{11}$ &   $E_1$, $E_2$, $E_3$, $E_4$, $E_6^*$   \\  \hline

& $[0 0 1 1 1 0 ]$ $[1 1 0 0 0 0]$ $[1 1 1 1 1 0]$  & $A_{3}$  &  $E_1$, $E_3$   \\  \hline

& $[0 0 1 1 1 1 ]$ $[1 1 0 0 1 0]$ $[1 1 1 0 1 0]$  & $A_{19}$  & $E_1$, $E_3$    \\  \hline

& $\bf [0 1 0 0 0 1]$ $[0 1 0 0 1 0]$ $[0 1 0 1 0 0]$ $[0 1 1 0 1 0]$ $[0 1 1 1 0 0]$ $\bf [0 1 1 1 1 1]$ $\bf [1 0 0 0 0 1]$ $[1 0 0 0 1 0 ]$ $[1 0 0 1 0 0]$ $[1 0 1 0 1 0]$ $[1 0 1 1 0 0]$ $\bf [1 0 1 1 1 1]$  &   $A_{22}$ &  $E_4$, $E_6^*$  \\  \hline

& $[0 1 0 0 1 1]$ $[0 1 0 1 0 1]$ $[0 1 0 1 1 0]$ $[0 1 1 0 0 0]$ $[0 1 1 0 1 1]$ $[0 1 1 1 0 1]$ $[1 0 0 0 1 1]$ $[1 0 0 1 0 1]$ $[1 0 0 1 1 0]$ $[1 0 1 0 0 0]$ $[1 0 1 0 1 1]$ $[1 0 1 1 0 1$ & $A_{23}$  &  $E_4$, $E_5$    \\  \hline

& $[0 1 0 1 1 1]$ $[0 1 1 0 0 1]$  $[1 0 0 1 1 1]$ $[1 0 1 0 0 1]$   & $A_{26}$  & $E_5$, $E_6^*$, $E_7$ \\  \hline

\end{tabular}
 \end{center}
\normalsize
 \begin{center} Table 4: Cosimple single-element coextensions of $D_1$, $D_2$, $D_3$ \end{center} 

\bigskip

%%%%%%%%%%%%%%%%%%%%%%%%%%%%%%%%%%%%%%%%%%%%%%%%%%%%%

\end {document}